\documentclass[11pt]{amsart}

\usepackage[english]{babel} 
\usepackage[T1]{fontenc}
\usepackage[utf8]{inputenc}

\usepackage[left=2.5cm, right=2.5cm, top=2cm, bottom=2cm]{geometry} 

\usepackage{amsmath,amssymb,mathrsfs} 
\usepackage{graphicx}  
\usepackage{float, caption} 
\usepackage{wrapfig} 

\usepackage{framed}
\usepackage{array, multirow, tabularx} 
\usepackage{tikz} 
\usepackage{pifont} 
\usetikzlibrary{matrix,arrows,decorations.pathmorphing,patterns}
\usepackage{mathtools}
\usepackage{enumerate}

\usepackage[bookmarksnumbered=true,pdfstartview=FitH,pdftex,pdfborder={0 0 0},linkcolor=blue,citecolor=blue,colorlinks=true]{hyperref} 

{\newtheorem{definition}{Definition}[section]}

{
\newtheorem{theorem}[definition]{\textsc{Theorem}}}
{
\newtheorem*{theoreme}{\textsc{Theorem}}}

{
\newtheorem{proposition}[definition]{Proposition}}

{
\newtheorem{lemma}[definition]{Lemma}}

{
\newtheorem{corollary}[definition]{Corollary}}

{\theoremstyle{remark}
\newtheorem{remark}[definition]{Remark}}

{\theoremstyle{remark}
\newtheorem{condition}[definition]{Conditions}}

{\theoremstyle{remark}
}

{\theoremstyle{remark}
\newtheorem{example}[definition]{Example}}

{\theoremstyle{remark}
\newtheorem{notation}[definition]{Notation}}

{
}

{

{}
{
}
\newcommand{\Hom}[3]{\mathrm{Hom}_{#1}\lrp{#2,#3}}
\newcommand{\extq}[4]{\mathrm{Ext}^{#1}_{#2}\lrp{#3,#4}}

\newcommand{\dr}{\ensuremath{\partial}}

\newcommand{\dd}{\ensuremath{\mathrm{d}}}

\newcommand{\eps}{\ensuremath{\varepsilon}}

\newcommand{\C}{\ensuremath{\mathbb{C}}}

\newcommand{\incl}{\ensuremath{\hookrightarrow}}

\newcommand{\ov}[1]{\ensuremath{\overline{#1}}}
\newcommand{\w}{\ensuremath{\omega}}

\newcommand{\co}{\ensuremath{\mathcal{O}}}

\newcommand{\wt}[1]{\ensuremath{\widetilde{#1}}}
\newcommand{\m}{\ensuremath{\mathfrak{m}}}

\newcommand{\f}[2]{\ensuremath{\frac{#1}{#2}}}

\newcommand{\fleche}{ edge node[auto]} 
\newcommand{\flecheb}{ edge node[auto, swap]}

\newcommand{\und}[1]{\underline{#1}}
\newcommand{\lra}[1]{\left\{#1\right\}}

\newcommand{\vphi}{\varphi}
\newcommand{\Z}{\mathbb{Z}}

\newcommand{\mc}[1]{\mathcal{#1}}

\newcommand{\derk}{\mathrm{Der}^k(-\log C)}

\newcommand{\homss}[2]{\mathrm{Hom}_{\co_S}\left(#1,#2\right)}

\newcommand{\homsi}[1]{\mathrm{Hom}_{\co_S}\left(#1,\Sigma\right)}
\newcommand{\wh}[1]{\widehat{#1}}

\newcommand{\resh}{\mathrm{res}_{C,\und{h}}}
\newcommand{\unk}{\{1,\ldots,k\}}
\newcommand{\res}{\mathrm{res}_C}
\newcommand{\der}{\mathrm{Der}}
\newcommand{\lrp}[1]{\left(#1\right)}

\newcommand{\extt}[3]{\mathrm{Ext}^{#1}_{\co_S}\left(#2,#3\right)}
\newcommand{\toto}[1]{\wt{\Omega}^{#1}}
\newcommand{\val}{\mathrm{val}}
\newcommand{\N}{\mathbb{N}}
\newcommand{\depth}{\mathrm{depth}}
\newcommand{\projdim}{\mathrm{projdim}}
\newcommand{\sm}{\smallskip}

\newcommand{\poin}{\bullet}
\newcommand{\ddh}{\mathrm{d}h_1\wedge \cdots \wedge \mathrm{d}h_k}
\newcommand{\ddf}{\mathrm{d}f_1\wedge \cdots \wedge \mathrm{d}f_k}
\newcommand{\omkxc}{\Omega^k(\log X/C)}
\newcommand{\totofq}[1]{\wt{\Omega}^{#1}_{\und{f}}}
\newcommand{\derkxc}{\mathrm{Der}^k(-\log X)}
\newcommand{\resxc}{\mathrm{res}_{X/C}}

\newcommand{\tomf}{\totofq{k}}

\newcommand{\lrb}[1]{\left\langle#1\right\rangle}

\usepackage{verbatim}

\title{Characterizations of freeness for equidimensional subspaces}
\author[D.~Pol]{Delphine Pol}
\email{\href{pol@mathematik.uni-kl.de}{pol@mathematik.uni-kl.de}}
\address{Department of Mathematics, T.U. Kaiserslautern\\
 Gottlieb-Daimler-Strasse building 48  \\
  67663 Kaiserslautern\\
  Germany }
 
 \date{\today}

\thanks{Research partially supported by a Japan Society for the Promotion of Science (JSPS) postdoctoral fellowship (Short-term) for North American and European researchers and by a postdoctoral fellowship from the Alexander Von Humboldt foundation.}

\begin{document}

\begin{abstract}
The purpose of this paper is to investigate properties of the minimal free resolution of the modules of multi-logarithmic forms along a reduced equidimensional subspace. We first consider a notion of freeness for reduced complete intersections, and more generally for reduced equidimensional subspaces embedded in a smooth manifold, which generalizes the notion of Saito free divisors.  The first main result is a characterization of freeness in terms of the projective dimension of the module of multi-logarithmic $k$-forms, where $k$ is the codimension. We also prove that there is a perfect pairing between the module of multi-logarithmic differential $k$-forms and the module of multi-logarithmic $k$-vector fields which generalizes the duality between the corresponding modules in the hypersurface case. We deduce from this perfect pairing a duality between the Jacobian ideal and the module of multi-residues of multi-logarithmic $k$-forms. In the last part of this paper, we investigate logarithmic modules along some examples of free singularities. The main result in this section is an explicit computation of the minimal free resolution of the module of multi-logarithmic forms and multi-residues for quasi-homogeneous complete intersection curves which uses our first main theorem.
\end{abstract}

\maketitle

\section{Introduction}
\label{intro}

Logarithmic differential forms along normal crossing divisors appear in P. Deligne's work, where he proves in particular that these modules are free, and in addition form a complex which computes the cohomology of the complement (see \cite{deligne-hodgeII}). The notion of logarithmic forms is then extended to arbitrary singular hypersurfaces by K. Saito in \cite{saito75} and \cite{saitolog}: a logarithmic differential form is a meromorphic form with simple poles along the hypersurface and such that its differential also has simple poles. He shows that the modules of logarithmic forms along the discriminant of a deformation of an isolated hypersurface singularity are free. However, this property is not always satisfied, and  the hypersurfaces for which the module of logarithmic $1$-forms is free are called \emph{free divisors}.

The question of identifying free divisors arises in several contexts,   for example in the study of prehomogeneous spaces (see \cite{granger-mond-schulze-prehomogeneous}), quivers (see \cite{buchweitz-mond-quiver}) or projective curves (\cite{dimca-sticlaru-nearly-free}, \cite{valles-pencil}). The study of freeness for hyperplane arrangements was initiated by H. Terao (see for example \cite{terao-hyperplanes-freeness-80}) and is still a very active field of research, in particular, Terao's conjecture asking if freeness only depends on the combinatorics is still open.

\medskip

In this paper we investigate a generalization of the notion of freeness for complete intersections, and more generally for equidimensional subspaces. 

\medskip

This generalization was first suggested in the case of complete intersections in \cite{gsci}, where the authors proposed to keep the term freeness since their definition is identical to  the characterization of free divisors by their singular locus proved by  H. Terao in \cite{terao-hyperplanes-freeness-80} and A.G. Aleksandrov in \cite{alek}. The definition of freeness suggested in \cite{gsci} is the following: a reduced complete intersection $C$ is called \textit{free} if it is smooth or if its singular locus is Cohen-Macaulay of codimension $1$ in $C$. 

\smallskip

\medskip

The following result reinforces the notion of freeness for equidimensional subspaces as a natural generalisation of freeness for hypersurfaces, since freeness means projective dimension $0$. Modules of multi-logarithmic forms are introduced in \cite{alekres} and \cite{aleksurvey} (see definition~\ref{multi-log}). We show in this paper that in all cases, freeness means that the module of multi-logarithmic $k$-form has minimal projective dimension, where $k$ is the codimension.  Our first main result is as follows (see also theorem~\ref{char:free}):

\begin{theoreme}
Let $X$ be a reduced equidimensional subspace of codimension $k$ in $(\C^m,0)$. 
\begin{itemize}
\item The projective dimension of the module of multi-logarithmic $k$-forms is greater than or equal to $k-1$.
\item In addition, the projective dimension is exactly $k-1$ if and only if $X$ is free. 
\end{itemize}
\end{theoreme}

The proof of this theorem uses duality properties proved in section \ref{sec:2}. In particular, we show the existence of a perfect pairing between the module of multi-logarithmic $k$-forms and the module of multi-logarithmic $k$-vector fields (see definition~\ref{def:derkxc}). This perfect pairing generalizes the duality between logarithmic $1$-forms and logarithmic vector fields along hypersurfaces proved in \cite{saitolog}. 

\smallskip

One of the fundamental example of free divisors is the normal crossing divisor. It is therefore natural to investigate freeness for the natural generalizations of normal crossing divisors, the so-called normal crossing singularities, which are an arbitrary equidimensional union of coordinate subspaces (see definition~\ref{de:normal:crossing}). We prove in this paper that normal crossing singularities of any codimension are free singularities. In particular, it means that the singular locus $\mathrm{Sing}(D)$ of the normal crossing divisor $D$ is free. Furthermore, from the point of view of hyperplane arrangements, it means that the subspace arrangement composed of the union of the subspaces of given rank in the intersection lattice of the Boolean arrangement is a free subspace arrangement. An interesting question would be to determine when such a property occurs.

\smallskip

Computing the modules of multi-logarithmic forms is in general a challenging question. Using the fact that reduced curves are free singularities,  we apply our main result to compute explicitly the minimal free resolution of the module of multi-logarithmic forms for reduced quasi-homogeneous complete intersection curves, which gives the second main result of this paper (see theorem~\ref{res:theo:omk:free:res}). We find that the Betti numbers in that case depend only on the codimension. Examples of free resolutions of the modules of multi-logarithmic forms of  homogeneous surfaces in $\C^4$ computed with \textsc{Singular} (\cite{DGPS}) show that the Betti numbers  depend on the surface (see example~\ref{ex:surf:c4}).

\smallskip

\sm

In addition, various homogeneous examples in $\C^4$, such as example~\ref{ex:surf:c4}, show that as in the case of hyperplane arrangements, it is a challenging question to investigate freeness of equidimensional central subspace arrangements, which are also related to projective arrangements.

\medskip

Let us describe now the content of this paper. 

\smallskip

We first give the basic definitions and properties of multi-logarithmic differential forms along reduced complete intersections and equidimensional subspaces. We also give an alternative proof to \cite[Proposition 2.1]{snc} of the inclusion of the ring of weakly holomorphic functions in the modules of logarithmic multi-residues, which is more elementary and does not use the isomorphism between multi-residues and regular meromorphic forms given in \cite{alektsikh}.

\smallskip 

In section~\ref{sec:nonci}, we consider a reduced equidimensional subspace $X$ of dimension $n$. Let $C$ be a reduced complete intersection of dimension $n$ defined by a regular sequence $(f_1,\ldots,f_k)$ containing $X$, and $f=f_1\cdots f_k$. We denote by $\mc{I}_C\subseteq \C\lra{x_1,\ldots,x_m}$ the ideal generated by $f_1,\ldots,f_k$. Let $c_X$ be the fundamental class of $X$ (see notation~\ref{nonci:nota:alpha0}). In particular, the fundamental class of the complete intersection $C$ is $\dd f_1\wedge \dots\wedge \dd f_k$. We  prove the following characterization of multi-logarithmic differential forms with respect to the pair $(X,C)$, which generalizes to equidimensional subspaces \cite[\S 3, theorem 1]{alekres} and \cite[(1.1)]{saitolog} (see proposition~\ref{nonci:prop:carac:loga}): a meromorphic $q$-form $\w\in\frac{1}{f}\Omega^q_S$ is multi-logarithmic if and only if there exist $g\in\co_S$ which induces a non zero divisor in $\co_C=\C\lra{x_1,\ldots,x_m}/\mc{I}_C$, $\xi\in\Omega^{q-k}_S$ and $\eta\in\frac{1}{f}\mc{I}_C\Omega^q$ such that $ g\w=\dfrac{c_X}{f}\wedge \xi+\eta$. 

\sm

We suggest the following definition of multi-logarithmic $k$-vector fields: a holomorphic $k$-vector field $\delta$ is multi-logarithmic if $\delta(c_X)\in\mc{I}_X$. We then have the following perfect pairing  which generalizes the duality of the hypersurface case (see proposition~\ref{nonci:prop:dualite:omkxc:derkxc}):
$$\derkxc\times\omkxc\to\frac{1}{f}\mc{I}_C.$$

In subsection~\ref{res:jac:sec}, we investigate duality properties of the module of multi-residues $\mc{R}_X$ of the multi-logarithmic $k$-forms. Let $\mc{J}_{X/C}$ denote the restriction to $X$ of the Jacobian ideal of $C$. Using an approach which is similar to the proof of \cite[Proposition 3.4]{gsres}, and which uses the previous perfect pairing, we prove that $\mathrm{Hom}_{\co_C}(\mc{J}_{X/C},\co_C)=\mc{R}_X$ (see proposition~\ref{res:jac}). In particular, for a complete intersection $C$, we have $\mathrm{Hom}_{\co_C}(\mc{J}_C,\co_C)=\mc{R}_C$, where $\mc{J}_C$ is the Jacobian ideal of $C$. 

\sm

Section~\ref{free:section} is devoted to the first main result of this paper, namely, the characterizations of freeness for an equidimensional subspace $X$ of codimension $k$. We first show the easy statement that the ideal $\mc{J}_{X/C}$ is Cohen-Macaulay if and only if the projective dimension of $\derkxc$ is $k-1$, which is a direct consequence of the Depth Lemma and the Auslander-Buchsbaum formula (see proposition~\ref{char:gs}). We then say that $X$ is free if $\mc{J}_{X/C}$ is Cohen-Macaulay, or equivalently, if the projective dimension of $\derkxc$ is $k-1$. This definition is compatible with the definition of freeness for complete intersections given in \cite{gsres}.
The modules of multi-logarithmic forms have as an homomorphic image the modules of multi-residues which are intrinsic and isomorphic to the modules of regular meromorphic forms. It is  therefore significant to look for a characterization of freeness involving the modules of multi-logarithmic forms.
Contrary to the hypersurface case, passing from $\derkxc$ to $\omkxc$ requires much more work.  The proof of our main theorem~\ref{char:free} is developed in section~\ref{free:section}, and uses Koszul complexes and change of rings spectral sequences. Our main theorem was recently reobtained by methods from pure commutative algebra introducing a residual duality over Gorenstein rings in \cite{schulze-tozzo-free}.

\sm

In section~\ref{normal-crossing-sec}, we show that the normal crossing singularities are free, by computing explicitely the modules of multi-logarithmic vector fields. 

\sm

As an application of theorem~\ref{char:free}, we obtain our second main result : we compute an explicit free resolution of the modules of multi-residues and multi-logarithmic forms for quasi-homogeneous complete intersection curves, which are free singularities  (see theorems~\ref{res:theo:res:res} and~\ref{res:theo:omk:free:res}). We use in this section results from \cite{polvalues}. This explicit computation is one of the key ingredients used in \cite{pol-solomon-terao} to study the case of line arrangements.

\medskip

\textbf{Acknowledgements}
The author is grateful to Michel Granger for raising several questions on this subject and many helpful discussions. The author also thanks Mathias Schulze for inviting her in Kaiserslautern in 2015, and for his suggestion to use Ischebeck's lemma instead of spectral sequences in section~\ref{res:jac:sec}, and for pointing out the characterization of quasi-homogenous curves used in the proof of proposition~\ref{gen:rc}. The author also thanks Takuro Abe for interesting discussions on freeness for subspace arrangements.

\section{Multi-logarithmic differential forms and multi-logarithmic vector fields}
\subsection{Complete intersections}

We will first consider the case of complete intersections, since the definitions for equidimensional subspaces rely on this case. 

Modules of multi-logarithmic differential forms along reduced complete intersections with arbitrary poles are introduced in \cite{alektsikhcrass}, and a variant with simple poles is introduced in \cite{alekres}. The precise relation between the two notions is described in \cite[\S 3.1.3]{polthese}. In order to work with finitely generated modules, we will consider the definition given in \cite{alekres}.

\begin{notation}
Let $m\geqslant 1$ be an integer and let $S=(\C^m,0)$. We denote by $\co_S$ the module of germs of holomorphic functions at the origin of $\C^m$. In particular, if $(x_1,\ldots,x_m)$ is a local system of coordinates of $S$, we identify $\co_S$ with $\C\lra{x_1,\ldots,x_m}$. For $q\in\N$, we denote by $\Omega^q$ the module of germs of holomorphic differential $q$-forms on $S$.
\end{notation}  

\begin{definition}[\cite{alekres}]
\label{multi-log}
Let $C\subseteq S$. We assume that the radical ideal $\mc{I}_C$ of vanishing functions on $C$ is generated by a regular sequence $(h_1,\ldots,h_k)\subseteq \co_S$, so that $C$ is a reduced complete intersection of codimension $k$ in $S$. We set $h=h_1\cdots h_k$. Let $q\in\N$. The module of \emph{multi-logarithmic $q$-forms on $C$ with respect to the equations $\underline{h}=(h_1,\ldots,h_k)$} is:
$$\Omega^q(\log C,\underline{h})=\lra{\w\in\frac{1}{h}\Omega^q\ ;\ \dd(\mc{I}_C)\wedge\w\subseteq \frac{1}{h}\mc{I}_C\Omega^{q+1}}.$$
\end{definition}

\begin{remark}
\label{rem:incl:log}
By \cite[Proposition 1]{alekres}, if $D$ is the hypersurface defined by $h$, for all $q\in\N$ we have the inclusion $\Omega^q(\log D)\subseteq \Omega^q(\log C,\und{h}).$
\end{remark}

\begin{remark}
\label{numerators}
From the definition, one can see that $h\cdot\Omega^q(\log C,\und{h})$ is the kernel of the map $\vphi:\Omega^q  \to \lrp{\dfrac{\Omega^{q+1}}{\mc{I}_C\Omega^{q+1}}}^k$ given by $\vphi(\w)=\lrp{\dd h_1\wedge \w,\ldots,\dd h_k \wedge \w}$. Therefore, the module of numerators $h\cdot\Omega^q(\log C,\und{h})$ does not depend on the choice of the equations $(h_1,\ldots,h_k)$. 
\end{remark}

The following characterization of multi-logarithmic forms generalizes \cite[(1.1)]{saitolog}:

\begin{theorem}[\protect{\cite[\S 3, theorem 1]{alekres}}]
\label{theo:alek}
Let $\w\in\dfrac{1}{h}\Omega^q_{S}$ with $q\in \N$. Then $\w$ is multi-logarithmic if and only if there exist a holomorphic function $g\in\co_{S}$ which  induces a non zero divisor in $\co_{C}=\co_S/\mc{I}_C$, a holomorphic differential form $\xi\in\Omega^{q-k}$ and a meromorphic $q$-form $\eta\in\frac{1}{h}\mc{I}_C\Omega^q$ such that:
\begin{equation}
\label{def:k:log}
g\w=\f{\dd h_1\wedge\dots\wedge \dd h_k}{h}\wedge\xi+\eta.
\end{equation}
\end{theorem}
\begin{remark}
\label{qtrop:petit}
For $q<k$, we have the equality $\Omega^q_S(\log C,\und{h})=\frac{1}{h}\mc{I}_C\Omega^q$. Using this property, one can compute a decomposition of the module of logarithmic $1$-forms of the hypersurface $D$ defined by $h$, which can be related to the splayedness condition in the case of codimension 2 as considered in  \cite{aluffifaber} (see \cite[Proposition 3.1.45, Corollary 3.1.49]{polthese} for more details).
\end{remark}

Theorem~\ref{theo:alek} enables us to define the modules of multi-residues as follows:

\begin{definition}[\protect{\cite[\S 4, Definition 1]{alekres}}]
\label{ci:de:res} Let  $\w\in\Omega^q(\log C,\und{h}), q\geqslant k$.  Let us assume that $g,\xi,\eta$ satisfy the properties of theorem~\ref{theo:alek}. Then the \emph{multi-residue} of $\w$ is: 
$$\resh(\w):={\f{\xi}{g}}\Big|_{C}\in \mc{M}_{C} \otimes_{\co_{C}} \Omega^{q-k}_{C}$$
where $\mc{M}_C$ is the sheaf of meromorphic functions on $C$ and for $p\in\N$, $\Omega^{p}_C={\f{\Omega^{p}}{\sum_{i=1}^k h_i\Omega^{p}+ \sum_{i=1}^k\dd h_i\wedge \Omega^{p-1}}}\Big|_C$ is the module of K\"ahler differential $p$-forms on $C$.
\end{definition}

The notion of multi-residue is well-defined with respect to the choices of $g$, $\xi$, $\eta$ in \eqref{def:k:log} (see \cite[\S 4 proposition 2]{alekres}). We denote $\mc{R}_{C,\und{h}}^{q-k}:=\resh\lrp{\Omega^q(\log C,\und{h})}$. If $q=k$, we set $\mc{R}_{C,\und{h}}:=\mc{R}_{C,\und{h}}^0\subset \mc{M}_C$. 

\begin{proposition}[\protect{\cite[\S 4, lemma 1]{alekres}}]
Let $q\in\N$. We have the following exact sequence of $\co_S$-modules:
\begin{equation}
\label{seq:res}
0\longrightarrow \frac{1}{h}\mc{I}_C\Omega^q \longrightarrow \Omega^q(\log C,\und{h}) \xrightarrow{\resh} \mc{R}^{q-k}_{C,\und{h}}\longrightarrow 0
\end{equation}
\end{proposition}

The following property gives the precise relation between the residue maps when we change the regular sequence defining $C$. In particular, it also gives a more elementary and precise proof of the fact that the modules of multi-residues do not depend on the choice of the defining equations which was proved in \cite{alekres} using the isomorphism between the modules of multi-residues and regular meromorphic forms given in \cite[Theorem 3.1]{alektsikh}.

\begin{proposition}
\label{res:eq}
Let  $(h_1,\ldots, h_k)$ and $(f_1,\ldots,f_k)$ be two regular sequences defining the same reduced germ of complete intersection $C$ of $S$. We set $f=f_1\cdots f_k$. Let $A=(a_{ij})_{1\leqslant i,j\leqslant k}$  be a $k\times k$ matrix with coefficients in $\co_S$ such that $(f_1,\ldots,f_k)^t=A(h_1,\ldots,h_k)^t$. Then for all $q\in\N$ and for all $\alpha\in f\cdot \Omega^q(\log C,\und{f})$:
$$\resh\lrp{\frac{\alpha}{h}}=\det(A) \mathrm{res}_{C,\und{f}}\lrp{\frac{\alpha}{f}}.$$

In particular, for all $p\geqslant 0$, the module $\mc{R}_{C,\und{h}}^p$ does not depend on the choice of the defining equations.
\end{proposition}
\begin{proof}
Let $\alpha\in f_1\cdots f_k\cdot\Omega^q(\log C,\und{f})$. Then there exists $g\in\co_S$ which induces a non zero divisor in $\co_C$, $\xi\in\Omega^{q-k}$ and $\eta\in\frac{1}{f}\mc{I}_C\Omega^q$ such that $g\alpha=\mathrm{d}f_1\wedge \dots \wedge \mathrm{d} f_k\wedge \xi+f \eta$. In addition, there exists $\nu\in \mc{I}_C\Omega^k$ such that: \begin{equation}
\label{ci:eq:ddf:ddh}
\ddf = \nu + \det(A) \ddh.
\end{equation}

 Thus, $g\alpha=\ddh\wedge \big(\det(A)\xi\big) +\nu\wedge\xi+f\eta$.
 
 Since $f\eta\in\mc{I}_C\Omega^q$ and $\nu\in\mc{I}_C\Omega^k$, we have $\displaystyle{\frac{\nu\wedge\xi+f\eta}{h_1\cdots h_k}}\in\frac{1}{h}\mc{I}_C\Omega^q$, so that $$\resh\lrp{\displaystyle{\frac{\alpha}{h}}}=\displaystyle{\frac{\det(A)\xi}{g}=\det(A)\mathrm{res}_{C,\und{f}}\lrp{\frac{\alpha}{f}}}.$$
 In addition, the conditions on $A$ implies that $\det(A)$ is a unit in $\co_S$ and in $\co_C$.
\end{proof}

\begin{notation}
In the rest of the paper, we will use the notation $\Omega^\bullet(\log C)$,  $\mc{R}_C^\bullet$ and $\res$ where the set of equations is implicitly $(h_1,\ldots,h_k)$. 
\end{notation}

\begin{notation}
\label{nota:jac}
We denote by $\mc{J}_{C}$ the Jacobian ideal of $C$, which is the ideal of $\co_{C}$ generated by the $k\times k$-minors of the Jacobian matrix of $(h_1,\ldots,h_k)$. The Jacobian ideal does not depend on the choice of the equations $(h_1,\ldots,h_k)$. 
\end{notation}

The following property can be deduced from \cite[Proposition 2.1]{snc} using the isomorphism of \cite[Theorem 3.1]{alektsikh}. We suggest here a more elementary proof which is similar to the proof of \cite[lemma 2.8]{saitolog} and does not require  \cite[Theorem 3.1]{alektsikh}.
\begin{proposition}[see \cite{gsci}]
\label{normaliz}
Let $\pi :\wt{C}\to C$ be the normalization of $C$.
We have the inclusion $\co_{\wt{C}}\subseteq \mc{R}_C$.
\end{proposition}
\begin{proof}  
Let $\alpha\in\co_{\wt{C}}$. Then, thanks to  \cite[Theorem 2]{lipman-sathaye}, for every $g\in\mc{J}_C$,  $g \alpha\in \co_C$. 

Thus, for every subset $J\subseteq\{1,\ldots, m\}$ with  $|J|=k$, there exists $a_J\in\co_S$ such that its class in $\co_C$ satisfies $\Delta_J \alpha=\ov{a_J}\in\co_C$, where $\Delta_J$ is the minor of the Jacobian matrix relative to the set $J$. Let $I,J$ be two subsets of $\{1,\ldots,m\}$ with $k$ elements. Then from the equality $\Delta_I\Delta_J \alpha-\Delta_J\Delta_I\alpha=0\in\co_C$ we deduce: 
$$\Delta_I a_J-\Delta_J a_I=h_1 b_1^{IJ}+\cdots + h_k b_k^{IJ}\in\co_S.$$

Let us define $\w=\f{\sum_{|J|=k} a_J \dd x_J}{h}\in\frac{1}{h}\Omega^k_S$. The previous equality gives:
\begin{align*}
\Delta_I\w&=\f{\sum_J\Delta_I a_J\dd x_J}{h}
          =a_I \f{\dd h_1\wedge\dots\wedge \dd h_k}{h} + \eta \text{ with } \eta\in\frac{1}{h}\mc{I}_C\Omega^k.
\end{align*}
Moreover, there exists a linear combination of the $\Delta_J$ which does not induce a zero-divisor in~$\co_C$ (see \cite{lipman-sathaye}). Therefore, $\w\in\Omega^k(\log C)$ and $\res(\w)=\alpha\in\mc{R}_{C}$.
\end{proof}

\subsection{Equidimensional subspaces}
\label{sec:nonci}

The modules of multi-logarithmic differential forms and multi-residues along  a reduced Cohen-Macaulay subspace are defined in \cite{aleksurvey}. The definitions are inspired by the description of regular meromorphic differential forms by residue symbols which has been introduced in \cite{kerskenregulare}, which is also used by M. Schulze in \cite{snc}.

\smallskip

Let $X\subseteq S$ be the germ of a reduced equidimensional analytic subset of $S$ of dimension $n$ defined by a radical ideal $\mc{I}_X$. We set $k=m-n$. 

One can prove that there exists a regular sequence $(f_1,\ldots,f_k)\subseteq \mc{I}_X$ such that the ideal $\mc{I}_C$ generated by $f_1,\ldots,f_k$ is radical (see \cite[Remark 4.3]{alektsikh} or \cite[Proposition 4.2.1]{polthese} for a detailed proof of this result). We fix such a sequence $(f_1,\ldots,f_k)$. We denote by $C$ the complete intersection defined by the ideal $\mc{I}_C$. In particular, $C=X\cup Y$, where $Y$ is of pure dimension $n$ and does not contain any component of $X$.  We set $f=f_1\cdots f_k$, and $\mc{I}_Y$ the radical ideal defining~$Y$.  

\begin{definition}[\protect{\cite[Definition 10.1]{aleksurvey}}]
Let $q\in\N$. We define the \emph{module of multi-logarithmic $q$-forms with respect to the pair $(X,C)$} as : 
$$\Omega^q(\log X/C)=\lra{\w\in\frac{1}{f_1\cdots f_k}\Omega^q\ ;\ \mc{I}_X\w \subseteq \frac{1}{f}\mc{I}_C\Omega^q \text{ and } (\dd \mc{I}_X)\wedge \w\subseteq\frac{1}{f}\mc{I}_C\Omega^{q+1}}.$$
\end{definition}

\begin{remark}
For all $q\in\N$, we have :
$$\frac{1}{f}\mc{I}_C\Omega^q\subseteq \Omega^q(\log X/C)\subseteq \Omega^q(\log C).$$

In addition, for a reduced complete intersection $C$, we have $\Omega^q(\log C/C)=\Omega^q(\log C)$.
\end{remark}

\begin{remark}
\label{remark:change:C}
Similarly to the complete intersection case, since $f\Omega^q(\log X/C)$ is the kernel of the map $ \beta :\Omega^q \to \lrp{\Omega^q/\mc{I}_C\Omega^q}^r\oplus \lrp{\Omega^{q+1}/\mc{I}_C\Omega^{q+1}}^r$
defined for $\w\in\Omega^q$ by $\beta(\w)=(h_1\w,\ldots,h_r\w, \dd h_1\wedge \w,\ldots ,\dd h_r \wedge \w)$ where $(h_1,\ldots,h_r)$ is a generating family of $\mc{I}_X$, the module of numerators $f\cdot \Omega^q(\log X/C)$ depends only on $X$ and $C$. The modules $f\cdot\Omega^q(\log X/C)$ depend on the choice of the complete intersection $C$, since the modules depend on the radical ideal defining $C$. 
\end{remark} 

\begin{definition}[\protect{\cite[Proposition 10.1]{aleksurvey}}]
The multi-residue map $\resxc : \Omega^q(\log X/C) \to \mc{M}_X\otimes_{\co_X} \Omega^{q-k}_X $ is defined as the restriction of the map $\res$ to $X$. 

For $q\in\N$, we set $\mc{R}_{X}^q:=\resxc(\Omega^{q+k}(\log X/C))$.
\end{definition}
\begin{remark} The modules of multi-residues along equidimensional subspaces are isomorphic to the modules of regular meromorphic forms (see  \cite[Theorem 10.2]{aleksurvey}). The proof uses \cite[(1.2)]{kerskenregulare}, \cite[(E.20)]{kunzkahler} and \cite[(E.21)]{kunzkahler}. In particular, it implies that for all $q\in\N$,  $\mc{R}_{X}^q$ does not depend on the choice of the complete intersection $C$ (see also \cite[Proposition 4.1.13]{polthese} for a detailed proof of this statement). 
\end{remark}
\begin{proposition}[\protect{\cite[Theorem 10.2]{aleksurvey}}]
\label{nonci:prop:suite:ex}
We have the following exact sequence of $\co_S$-modules:
\begin{equation}
\label{nonci:eq:suite:ex}
0\to \frac{1}{f}\mc{I}_C\Omega^q\to\Omega^q(\log X/C)\xrightarrow{\resxc} \mc{R}_{X}^{q-k}\to 0.
\end{equation}
\end{proposition}

We prove in the following proposition~\ref{nonci:prop:carac:loga} a characterization of multi-logarithmic $q$-forms with respect to the pair $(X,C)$, which generalizes theorem~\ref{theo:alek}. We first need to introduce the notion of  fundamental form ${c_X}$ of $X$ (see for example \cite[(1.3)]{kerskenregulare}). We recall that $C=X\cup Y$ is an irredundant decomposition of $C$. 
\begin{notation}
\label{beta}
Let $\beta_f\in\co_{\wt{C}}$ be such that ${\beta_f}|_{X}=1$ and ${\beta_f}|_{Y}=0$. 
\end{notation}

A direct consequence of \cite{lipman-sathaye} is that the form $\ov{\beta_f \dd f_1\wedge\dots \wedge \dd f_k}\in\Omega^k_S\otimes_{\co_S}\mc{M}_C$ satisfies $\ov{\beta_f\dd f_1\wedge \dots \wedge \dd f_k}\in\Omega^k_S\otimes_{\co_S}\co_C$.

\begin{notation}
\label{nonci:nota:alpha0}
A fundamental form of $X$ is a form $c_{X}\in\Omega^k$ such that $\ov{c_{X}}=\ov{\beta_f\dd f_1\wedge \dots\wedge \dd f_k}\in\Omega^k\otimes_{\co_S}\co_C$. In particular, $\dd f_1\wedge \dots \wedge \dd f_k$ is a fundamental form of the complete intersection $C$. 
\end{notation}
\begin{remark}
The choice of the form $c_X$ is not unique, and depends also on the choice of the complete intersection $C$. However, all the possible choice of fundamental forms of $X$ as in the previous notation give the same element in the ring of regular meromorphic $0$-form $\w_X^0$ (see \cite{kerskenregulare}). 
\end{remark}

We fix a fundamental form $c_X$ of $X$. The following result is a consequence of the definition of $c_X$ and \cite[(1.3)]{kerskenregulare}: 
\begin{proposition}[see \protect{\cite[(1.3)]{kerskenregulare}}]
We have $\frac{c_X}{f}\in\Omega^k(\log X/C)$ and $\resxc\lrp{\frac{c_X}{f}}=1\in\mc{M}_X$. 
\end{proposition}

\begin{proposition}
\label{nonci:prop:carac:loga}
Let $\w\in\frac{1}{f}\Omega^q_S$. Then $\w\in\Omega^q(\log X/C)$ if and only if there exist $g\in\co_S$ which induces a non zero divisor in $\co_C$, $\xi\in\Omega^{q-k}_S$ and $\eta\in\frac{1}{f}\mc{I}_C\Omega^q$ such that: 
\begin{equation}
\label{nonci:eq:carac:loga}
g\w=\frac{c_X}{f}\wedge \xi+\eta
\end{equation}

We then have $\resxc\lrp{\w}={\frac{\xi}{g}}\big|_{X}$. 
\end{proposition}
\begin{proof}
Let $\w\in\Omega^q(\log X/C)$. Then $\w\in\Omega^q(\log C)$. Let $g,\xi,\eta$ satisfying theorem~\ref{theo:alek} such that: 
$$g\w=\frac{\dd f_1\wedge\dots \wedge \dd f_k}{f}\wedge \xi+\eta.$$
Then $g\resxc(\w)=\xi=\resxc\lrp{\frac{c_X}{f}\wedge\xi}$. By proposition~\ref{nonci:prop:suite:ex}, there exists $\eta'\in\frac{1}{f}\mc{I}_C\Omega^q$ such that $$g\w=\frac{c_X}{f}\wedge \xi +\eta'.$$

Conversely, let $\w\in\frac{1}{f}\Omega^q$ be such that $g\w=\frac{c_X}{f}\wedge \xi+\eta$ with $g,\xi,\eta$ as in the statement of the proposition.  Let $h\in\mc{I}_X$. Since $\frac{c_X}{f}\in\Omega^k(\log X/C)$, we deduce that $hg\w\in\frac{1}{f}\mc{I}_C\Omega^q$ and $\dd h\wedge g\w\in\frac{1}{f}\mc{I}_C\Omega^{q+1}$.  Therefore, since $g$ induces a non zero divisor in $\co_C$, we have $\w\in\Omega^q(\log X/C)$. 
\end{proof}

We end this section with the following lemma, which enables us to identify elements in $\Omega^q(\log C)$ which belongs to $\Omega^q(\log X/C)$, and which is quite useful for computations with \textsc{Singular}.
\begin{lemma}
\label{omegaC:IY}
For all $q\in\N$ we have: 
$$\Omega^q(\log X/C)= \Omega^q(\log C)\cap \frac{1}{f}\mc{I}_Y\Omega^q.$$
\end{lemma}
\begin{proof}
Let us prove the inclusion $\subseteq$. We already mention the inclusion $\Omega^q(\log X/C)\subseteq \Omega^q(\log C)$. Let $\w\in\Omega^q(\log X/C)$. Since $\mc{I}_X\w\subseteq \frac{1}{f}\mc{I}_C\Omega^q$, and $C=X\cup Y$ is an irredundant decomposition of $C$, we have $\w\in\frac{1}{f}\mc{I}_Y\Omega^q$. Let us prove the converse inclusion. Let $\w\in\Omega^q(\log C)\cap \frac{1}{f}\mc{I}_Y\Omega^q$. Then since $C=X\cup Y$, $\mc{I}_X\w\subseteq \frac{1}{f}\mc{I}_C\Omega^q$. Let $h\in\mc{I}_X$ and $F\in\mc{I}_Y$ be such that $F$ induces a non zero divisor in $\co_X$. Then $hF\in\mc{I}_C$, so that $\dd (hF)\wedge\w=\dd h\wedge F\w+\dd F\wedge h\w\in\frac{1}{f}\mc{I}_C\Omega^q$. Since $\mc{I}_C\subseteq \mc{I}_X$ and $h\in\mc{I}_X$, we have $\dd h \wedge F\w\in\frac{1}{f}\mc{I}_X\Omega^q$, and since $F$ is a non zero divisor in $\co_X$, we have $\dd h\wedge \w\in\frac{1}{f}\mc{I}_X\Omega^q$. Since $\w\in\frac{1}{f}\mc{I}_Y\Omega^q$ and $C=X\cup Y$, we have  $\dd h\wedge\w\in\frac{1}{f}\mc{I}_C\Omega^q$.
\end{proof}

\begin{remark}
\label{remark:residus:nul:Y}
Let $\w\in\Omega^q(\log X/C)$ and $g, \xi, \eta$ as in theorem~\ref{theo:alek} such that $g\w=\frac{\dd f_1\wedge\dots\wedge\dd f_k}{f}\wedge \xi+\eta$. One can deduce from lemma~\ref{omegaC:IY} that $\xi\in \mc{I}_Y$. 
\end{remark}

\subsection{Multi-logarithmic vector fields}

We introduce here a module of multi-logarithmic $k$-vector fields. We will prove in proposition~\ref{nonci:prop:dualite:omkxc:derkxc} that there exists a perfect pairing between this module and $\omkxc$  with values in $\frac{1}{f}\mc{I}_C$, which generalizes \cite[(1.6)]{saitolog}.

Let $\Theta_S$ denote the $\co_S$-module of holomorphic vector fields and $\Theta^k_S:=\bigwedge^k \Theta_S$ be the exterior power of order $k$ of $\Theta_S$. 

We can evaluate  a $k$-form $\w\in\Omega^k_S$ on a $k$-vector field $\delta\in\Theta^k_S$, which gives a function denoted by $\w(\delta)$ or $\delta(\w)$.

We suggest the following definition of multi-logarithmic $k$-vector fields, which generalizes  \cite[(5.1)]{gsci}:
\begin{definition}
\label{def:derkxc}
Let $\delta\in\Theta^k_S$. We say that $\delta$ is a \emph{multi-logarithmic $k$-vector field with respect to $X$} if $c_X(\delta)\in\mc{I}_X$. We denote by $\derkxc$ the module of multi-logarithmic $k$-vector fields along~$X$. 
\end{definition}
\begin{remark}
An equivalent definition of $\derkxc$ involving only $\mc{I}_X$ is given in \cite[Definition 3.7]{schulze-tozzo-free}, which shows that $\derkxc$ is independent from the choice of $C$, property which can also be proved directly with our definition (see remark~\ref{remark:change:ic:jxc}). 
\end{remark}

 \begin{remark}
 Since $c_X\in\mc{I}_Y\Omega^k$, for all $\delta\in\Theta^k_S$, $c_X(\delta)\in\mc{I}_Y$. Therefore, $$\derkxc=\lra{\delta\in\Theta^k_S\ ;\ c_X(\delta)\in\mc{I}_C}.$$
 \end{remark}
 
\begin{remark}
If $C$ is a reduced complete intersection, this definition coincides with \cite[(5.1)]{gsci}.
\end{remark} 
 As a consequence, we have the following proposition:
 \begin{proposition}
 We have the following inclusion: 
 $$\derk\subseteq\derkxc.$$
 \end{proposition}

\begin{remark}
The notion of logarithmic vector field studied in \cite{saitolog} can be extended to spaces of higher codimension (see \cite{hauser-muller}).  A holomorphic vector field $\eta\in\Theta_S$ is logarithmic along an equidimensional space $X$ if $\eta$ is tangent to $X$ at its smooth points. Even in the complete intersection case, we cannot obtain in general all the multi-logarithmic $k$-vector fields from the logarithmic vector fields. Let us consider a reduced complete intersection $C$ defined by the ideal $\mc{I}_C$ generated by the regular sequence $(h_1,\ldots,h_k)$. The module of logarithmic vector fields is $\der(-\log C)=\lra{\eta\in\Theta_S\ ;\ \eta(\mc{I}_C)\subseteq \mc{I}_C}$. The following inclusion is a direct consequence of the definition:
\begin{equation}
\label{ci:eq:der}
\der(-\log C)\wedge  \Theta_S^{k-1}\subseteq \derk.
\end{equation}

By considering for example the complete intersection curve defined by $h_1=x^3-y^2$ and $h_2=x^2y-z^2$ (see \cite[Exemple 3.2.10]{polthese}), one can check that the inclusion~\eqref{ci:eq:der} may be strict, which answers a question of Luis Narv\'aez-Macarro.
\end{remark}

\section{Duality results}
\label{sec:2}

The proof of our main theorem~\ref{char:free} requires the following duality results, which are generalizations to complete intersections and equidimensional subspaces of the dualities proved in \cite[Lemma (1.6)]{saitolog} and \cite[Proposition 3.4]{gsres} for hypersurfaces. 

\subsection{Perfect pairing between multi-logarithmic vector fields and multi-logarithmic forms}
For a reduced hypersurface $D$, a $\co_S$-duality is satisfied between $\Omega^1(\log D)$ and $\der(-\log D)$ (see \cite[Lemma (1.6)]{saitolog}). 

Let us prove that there is a perfect pairing between $\omkxc$ and $\derkxc$ which generalizes the duality of the hypersurface case, and whose proof is inspired by the proof of \cite[Lemma (1.6)]{saitolog}.

\begin{notation}
To simplify the notations, we set for $q\in\N$, $\totofq{q}=\frac{1}{f}\mc{I}_C\Omega^q$ and we set $\Sigma=\frac{1}{f}\mc{I}_C=\totofq{0}$.
\end{notation}

We first need the following lemma, which is a direct consequence of proposition~\ref{nonci:prop:carac:loga} and the definition of $\derkxc$. 
\begin{lemma}
\label{nonci:lem:pairing}
Let $\delta\in\derkxc$ and $\w\in\omkxc$. Then $\w(\delta)\in\Sigma$. 
\end{lemma}

Thanks to lemma~\ref{nonci:lem:pairing}, we see that we have a natural pairing $$\derkxc\times \omkxc\to \Sigma.$$

The following lemma is used in the proof of proposition~\ref{nonci:prop:dualite:omkxc:derkxc}:
\begin{lemma}
\label{ci:lem:dual:toto}
We have the following perfect pairings: 
\begin{align}
\totofq{k}\times\Theta^k_S&\to \Sigma,\\
\frac{1}{f}\Omega^k\times \sum_{i=1}^k f_i \Theta^k_S&\to \Sigma.
\end{align}
\end{lemma}
\begin{proof}
Let us notice that $\totofq{k}=\Omega^k\otimes \Sigma$ and $\sum f_i\Theta^k_S=f\Theta^k_S\otimes \Sigma$. Let $M$ be either $\Omega^k_S$ or $f\Theta^k_S$. It is easy to prove that $\mathrm{Hom}_{\co_S}\lrp{M\otimes\Sigma,\Sigma}=\mathrm{Hom}_{\co_S}\lrp{M,\co_S}$ and $\mathrm{Hom}_{\co_S}\lrp{M,\Sigma}=\mathrm{Hom}_{\co_S}\lrp{M,\co_S}\otimes\Sigma$. Hence the result.
\end{proof}
\begin{proposition}
\label{nonci:prop:dualite:omkxc:derkxc}
We have the following perfect pairing: 
\begin{equation}
\label{nonci:eq:perfect:pairing:omkxc:derkxc}
\Omega^k(\log X/C)\times \derkxc\to \Sigma.
\end{equation}
\end{proposition}
\begin{proof}
We have the following inclusions:
\begin{equation}
\label{nonci:eq:suite:dualite:om}
\totofq{k}\subseteq \omkxc\subseteq \frac{1}{f}\Omega^k_S
\end{equation}
\begin{equation}
\label{nonci:eq:suite:dualite:der}
\Theta^k_S\supseteq \derkxc\supseteq \mc{I}_C\Theta^k_S.
\end{equation}
We deduce from lemmas~\ref{nonci:lem:pairing}, \ref{ci:lem:dual:toto} and the inclusions~\eqref{nonci:eq:suite:dualite:om} and~\eqref{nonci:eq:suite:dualite:der} the following: 
\begin{align*}
\derkxc\subseteq\mathrm{Hom}_{\co_S}\lrp{\omkxc,\Sigma}\subseteq \Theta^k_S,\\
\omkxc\subseteq\mathrm{Hom}_{\co_S}\lrp{\derkxc,\Sigma}\subseteq \frac{1}{f}\Omega^k_S.
\end{align*}
Let us prove that the left-hand-side inclusions are equalities. 

Since $\frac{c_X}{f}\in\omkxc$, for all $\delta\in \mathrm{Hom}_{\co_S}\lrp{\omkxc,\Sigma}$ we have $c_X(\delta)\in\mc{I}_C$ so that $\delta\in\derkxc$. 
 Therefore, $\mathrm{Hom}_{\co_S}\lrp{\omkxc,\Sigma}=\derkxc$. 
 
 Let $\w\in\mathrm{Hom}_{\co_S}\lrp{\derkxc,\Sigma}$. Let us prove that $\w\in\omkxc$. Let us set $\w=\sum_{|I|=k} \frac{1}{f}\w_I\dd x_I$. Let $h\in\mc{I}_X$. Then for all $I$, $h\dr_{x_I}\in\derkxc$, so that $h\dr_{x_I}(\w)=h\frac{1}{f}\w_I\in\Sigma$. Therefore, $h\w\in\totofq{k}$. 

We set for $J\subseteq \lra{1,\ldots,m}$ with $|J|=k+1$,
$\delta_J=\sum_{\ell=1}^{k+1} (-1)^{\ell-1} \frac{\dr h}{\dr x_{j_\ell}} \dr x_{j_1}\wedge \dots\wedge\wh{\dr x_{j_\ell}}\wedge \dots\wedge\dr x_{j_{k+1}}$.

We thus have $\w(\delta_J)=(\dd h\wedge\w)(\dr_{x_J})$. Let us prove that $\delta_J\in\derkxc$. We have: 
$$\lrp{\beta_f\dd f_1\wedge\dots\wedge \dd f_k}(\delta_J)=\beta_f\sum_{\ell=1}^{k+1} (-1)^{\ell-1} \f{\dr h}{\dr x_{j_\ell}} \Delta_{j_1\dots \widehat{j_\ell}\dots j_{k+1}}=\beta_f\begin{vmatrix}
\frac{\dr h}{\dr x_{j_1}} & \dots & \frac{\dr h}{\dr x_{j_{k+1}}}\\
\frac{\dr f_1}{\dr x_{j_1}} & \dots & \frac{\dr f_1}{\dr x_{j_{k+1}}}\\
\vdots & & \vdots\\
\frac{\dr f_k}{\dr x_{j_1}}& \dots & \frac{\dr f_k}{\dr x_{j_{k+1}}}
\end{vmatrix}.$$

Since the codimension of $X$ is $k$ and $(h,f_1,\ldots,f_k)\subseteq \mc{I}_X$, the restriction of the previous determinant to $X$ is zero. Given that ${\beta_f}|_{Y}=0$, we have $\ov{(\beta_f\dd f_1\wedge\dots\wedge \dd f_k)(\delta_J)}=\ov{0}\in\mc{M}_C$. Thus, $\delta_J\in\derkxc$. Since $\w(\delta_J)=(\dd h\wedge\w)(\dr_{x_J})$, we deduce that $\dd h\wedge\w\in\totofq{k+1}$ and $\mathrm{Hom}_{\co_S}\lrp{\derkxc,\Sigma}=\omkxc$.
\end{proof}

\begin{remark}
By lemma~\ref{ci:lem:dual:toto} and proposition~\ref{nonci:prop:dualite:omkxc:derkxc}, if $C$ is a reduced complete intersection, we have the following perfect pairings between modules which all depend only on~$C$ and not on the choice of equations: 
\begin{align*}
\mc{I}_C\Omega^k\times \Theta^k_S&\to \mc{I}_C\\
\Omega^k\times \mc{I}_C\Theta^k_S&\to \mc{I}_C\\
f\cdot \Omega^k(\log C)\times \derk&\to \mc{I}_C
\end{align*}
\end{remark}

\subsection{Multi-residues and Jacobian ideal}
\label{res:jac:sec}

Let $C$ be a reduced complete intersection. A consequence of \cite[Lemma 5.4]{snc} and \cite[Theorem 3.1]{alektsikh} is that the dual of the Jacobian ideal $\mc{J}_C$ of $C$ is the module of multi-residues $\mc{R}_C$. 

In this subsection, we give another proof of this duality, which does not depend on the isomorphism of \cite[Theorem 3.1]{alektsikh}. Furthermore, our approach enables us to extend this duality to the case of reduced equidimensional subspaces by introducing the ideal $\mc{J}_{X/C}$ defined below which may differ from the ideal considered in \cite{snc} when the subspace is not Gorenstein (see remark~\ref{remar:jac}). Our proof uses the perfect pairings of proposition~\ref{nonci:prop:dualite:omkxc:derkxc} and is analogous to the proof of \cite[Proposition 3.4]{gsres}. 

\begin{notation}
Let $X$ be a reduced equidimensional subspace of codimension $k$ in $S$ and $C$ be a reduced complete intersection of codimension $k$ containing $X$.  We denote by $\mc{J}_{X/C}\subseteq \co_X$ the restriction of the Jacobian ideal $\mc{J}_C$ of $C$ to the space $X$.
\end{notation}
\begin{remark}
Since ${\beta_f}|_{X}=1$, the ideal $\mc{J}_{X/C}$ is given by $\mc{J}_{X/C}=\lra{\ov{\delta(c_X)}\in\co_X\ ;\ \delta\in\Theta^k_S}$. 
\end{remark}
\begin{remark}
\label{remark:change:ic:jxc}
If $C''$ is another reduced complete intersection of codimension $k$ containing $X$, the ideals $\mc{J}_{X/C}$ and $\mc{J}_{X/C''}$ are isomorphic. To see this, let us consider a reduced complete intersection $C'$ containing $C$ and $C''$. Let  $(f_1',\ldots,f_k')$ be a regular sequence defining $C'$ and $A$ be a transition matrix from $(f_1',\ldots,f_k')$ to $(f_1,\ldots,f_k)$. There exists $\nu\in\mc{I}_C\Omega^k$ such that:
$$\dd f_1'\wedge\dots\wedge\dd f_k'=\det(A)\dd f_1\wedge\dots\wedge\dd f_k+\nu.$$

Therefore, with the notations of \ref{beta}, the image of $\beta_{f'}\dd f_1'\wedge\dots\wedge \dd f_k'$ in $\Omega^k\otimes_{\co_S}\co_X$ is:
$$\ov{\beta_{f'}\dd f_1'\wedge \dd f_k'}=\ov{\beta_{f'}\dd f_1'\wedge \dots\wedge \dd f_k'}=\ov{\beta_{f'}\det(A)\dd f_1\wedge \dots \wedge \dd f_k.}$$

As in notation~\ref{nonci:nota:alpha0}, let $c_{X}'\in\Omega^k$ be such that $\ov{c_X'}=\ov{\beta_{f'}\dd f_1'\wedge \dd f_k'}\in\Omega^k\otimes\co_{C'}$.
Since $\beta_{f'}|_X=\beta_{f}|_X=1$, we have $\ov{\delta(c_X')}=\det(A)\ov{\delta(c_X)}\in\co_X$, so that $\mc{J}_{X/C'}=\det(A)\ov{\mc{J}_{X/C}}$. In addition, since $X\subseteq C\cap C'$ and $X$ is not included in the singular locus of $C$ or $C'$, $\det(A)$ is a non zero divisor of $\co_X$, so that $\mc{J}_{X/C'}$ and $\mc{J}_{X/C}$ are isomorphic. It also shows directly with our definition~\ref{def:derkxc} that the module $\derkxc$ does not depend on the choice of $C$. 
\end{remark}

The following result is a consequence of the definitions:
\begin{proposition}
We have the following exact sequence of $\co_S$-modules:
$$0\to\derkxc\to\Theta^k_S\to \mc{J}_{X/C}\to 0$$
 where the last map is given by $\delta\in\Theta^k\mapsto \delta(c_X)\in\co_X$.  
\end{proposition}

\begin{remark}
\label{remar:jac}
Several notions of Jacobian ideals are associated with a reduced equidimensional subspace: the ideal $\mc{J}_X\subseteq \co_X$ generated by the $k\times k$ minors of the Jacobian matrix of $(h_1,\ldots,h_r)$, where $\sum_{i=1}^r h_i\co_S=\mc{I}_X$, the $\w$-Jacobian $\mc{J}_X^\w$ (see for example \cite{snc}), and the ideal $\mc{J}_{X/C}$ considered above. For the non Gorenstein curve of $\C^3$ defined by $h_1=y^3-x^2z$, $h_2=x^3y-z^2$ and $h_3=x^5-y^2z$, which is the irreducible curve parametrized by $(t^5,t^7,t^{11})$, one can check that $\mc{J}_X$, $\mc{J}_{X/C}$ and $\mc{J}_X^\w$ are pairwise distinct (see \cite[Exemple 4.2.35]{polthese}). 

\end{remark}

 We first need the following lemma:

\begin{lemma} 
\label{ext:jac}
Let $X$ be a reduced space of pure dimension $n=m-k$. We assume $k\geqslant 2$.  Then $\mathrm{Ext}^1_{\co_S}\lrp{\mc{J}_{X/C},\co_S}=0$ and $\mathrm{Ext}^1_{\co_S}\lrp{\mc{J}_{X/C},\Sigma}=\mathrm{Hom}_{\co_C}\lrp{\mc{J}_{X/C},\co_C}$.
\end{lemma}
\begin{proof}
We apply the functor $\mathrm{Hom}_{\co_S}(\mc{J}_{X/C},-)$ to the exact sequence $0\to\Sigma\xrightarrow{\times f} \co_S \to \co_C\to 0$.

It gives: $$0\to\mathrm{Hom}_{\co_S}(\mc{J}_{X/C},\co_C)\to \mathrm{Ext}^1_{\co_S}\lrp{\mc{J}_{X/C},\Sigma}\to \mathrm{Ext}^1_{\co_S}\lrp{\mc{J}_{X/C},\co_S}\to\dots$$

The depth of $\co_S$ is $m$ and since $\mc{J}_{X/C}$ is a fractional ideal of $\co_X$, the dimension of $\mc{J}_{X/C}$ is $m-k=\dim \co_X$. Thus, by Ischebeck's lemma (see \cite[15.E]{matsumura}), we have $\mathrm{Ext}_{\co_S}^1\lrp{\mc{J}_{X/C},\co_S}=0$. Hence the result.
\end{proof}

\begin{notation}
Let $I\subset \mc{M}_C$ be an ideal. We set $I^\vee=\mathrm{Hom}_{\co_C}(I,\co_C)$. 
\end{notation}

\begin{proposition}
\label{res:jac}
We have $\mc{J}_{X/C}^\vee\simeq\mc{R}_X$. In particular, if $C$ is a reduced complete intersection, $\mc{J}_C^\vee=\mc{R}_C$. 
\end{proposition}
\begin{proof}
We assume $k\geqslant 2$. For $k=1$, we refer to \cite[Proposition 3.4]{gsres}. 

We consider the double complex $\mathrm{Hom}_{\co_S}\lrp{\derkxc\incl\Theta^k_S, f : \Sigma\to\co_S}$, which gives almost the same diagram as the dual of $(3.8)$ in~\cite{gsres}. By lemma~\ref{ext:jac}, $\mathrm{Ext}^1_{\co_S}(\mc{J}_{X/C},\co_S)=0$, so that we obtain the following commutative diagram:

{\small
\begin{center}
\begin{tikzpicture}[scale=0.1][description/.style={fill=white,inner sep=1pt}]
\matrix (m) [matrix of math nodes, row sep=1.5em,
column sep=0.4em, text height=1.5ex, text depth=0.25ex]
{ & & & & 0   & & 0   & &   & &  \\
& & 0  & & \homss{\Theta^k_S}{\Sigma}    & & \homsi{\derkxc}    & & \extt{1}{\mc{J}_{X/C}}{\Sigma}   & & 0  \\
& & 0 & & \homss{\Theta^k_S}{\co_S}    & &  \mathrm{Hom}_{\co_S}\lrp{\derkxc,\co_S}   & & 0   & &    \\
0 & & \homss{\mc{J}_{X/C}}{\co_C}  & &  \homss{\Theta^k_S}{\co_C}    & & \homss{\derkxc}{\co_C}   & &  \extt{1}{\mc{J}_{X/C}}{\co_C} & & 0\\
& &\extt{1}{\mc{J}_{X/C}}{\Sigma}  & & 0    & &  \extt{1}{\derkxc}{\Sigma}   & & \extt{2}{\mc{J}_{X/C}}{\Sigma} & & 0 \\
 & & 0 && && && && \\};
\path[->,font=\scriptsize]
(m-2-3) \fleche {} (m-2-5)
(m-2-5) \fleche {} (m-2-7)
(m-2-7) \fleche {} (m-2-9)
(m-2-9) \fleche {} (m-2-11)
(m-3-3) \fleche {} (m-3-5)
(m-3-5) \fleche {} (m-3-7)
(m-3-7) \fleche {} (m-3-9)
(m-4-1) \fleche {} (m-4-3)
(m-4-3) \fleche {} (m-4-5)
(m-4-5) \fleche {} (m-4-7)
(m-4-7) \fleche {} (m-4-9)
(m-4-9) \fleche {} (m-4-11)
(m-5-3) \fleche {} (m-5-5)
(m-5-5) \fleche {} (m-5-7)
(m-5-7) \fleche {} (m-5-9)
(m-5-9) \fleche {} (m-5-11)
(m-1-5) \fleche {} (m-2-5)
(m-1-7) \fleche {} (m-2-7)
(m-2-5) \fleche {} (m-3-5)
(m-2-7) \fleche {$\cdot f$} (m-3-7)
(m-2-9) \fleche {} (m-3-9)
(m-3-5) \fleche {} (m-4-5)
(m-3-7) \fleche {} (m-4-7)
(m-3-9) \fleche {} (m-4-9)
(m-4-5) \fleche {} (m-5-5)
(m-4-7) \fleche {} (m-5-7)
(m-4-9) \fleche {} (m-5-9)
(m-4-3) \fleche {} (m-5-3)
(m-3-3) \fleche {} (m-4-3)
(m-5-3) \fleche {} (m-6-3);
\end{tikzpicture}
\end{center}}

Let $\vphi : \derkxc\to\Sigma$. By identifying $\mathrm{Hom}_{\co_S}\lrp{\derkxc,\Sigma}$ with $\omkxc$ thanks to proposition~\ref{nonci:prop:dualite:omkxc:derkxc},  we associate with $\vphi$ the form $\w=\f{1}{f}\sum_I \vphi(f\dr x_I)\dd x_I\in\omkxc$. We have $\vphi(\delta)=\w(\delta)$.

By a diagram chasing process, we obtain the map:
$$\begin{aligned} \omkxc & \to \Hom{\co_S}{\mc{J}_{X/C}}{\co_C}\\
 \w&\mapsto \lrp{\ov{a}\mapsto \ov{\res(\w) a}}.
 \end{aligned}$$

The map $\ov{a}\mapsto \ov{\res(\w)a}\in\co_C$ is well defined since by remark~\ref{remark:residus:nul:Y}, $\res(\w)|_Y=0$. 

Similarly, the same diagram chasing process starting from the lower left $\homss{\mc{J}_{X/C}}{\co_C}$ to the  upper right $\homss{\derkxc}{\Sigma}$, show that the map 
\begin{align*}
\theta : \mc{R}_X &\to \mc{J}_{X/C}^\vee\\
                   \rho & \mapsto \theta_\rho : \left\{\begin{aligned}
                                                                                \mc{J}_{X/C} & \to  \co_C\\
                                                                                    \ov{a} & \mapsto \ov{\rho a}
                                                                                    \end{aligned}\right.                                                                                 
\end{align*}
is an isomorphism.
\end{proof}

\section{Freeness for equidimensional subspaces}
\label{sec:3}
\label{free:section}

We prove here our main result, namely, theorem~\ref{char:free}, which is a characterization of freeness for equidimensional subspaces by the minimality of the projective dimension of the module $\omkxc$, which generalizes the hypersurface case. 

\subsection{Definition and statements}

The purpose of this section is to develop an analogue of freeness for equidimensional subspaces, and in particular for complete intersection, which generalizes the notion of Saito free divisors. A hypersurface is called free if the module of logarithmic vector fields is free. Among the different characterizations of free divisors, let us mention the following one:

\begin{theorem}[\protect{\cite{terao-hyperplanes-freeness-80}}, \protect{\cite{alek}}]
\label{char:alek}
The germ of a reduced singular divisor is free if and only if its singular locus is Cohen Macaulay of codimension $1$ in $D$. Equivalently, a reduced divisor $D$ is free if and only if the Jacobian ideal of $D$ is Cohen-Macaulay.
\end{theorem}

 Let us first notice the following property:

\begin{proposition}
\label{char:gs}
Let $X$ be a reduced equidimensional subspace of codimension $k$ in $S$, and $C$ be a reduced complete intersection of codimension $k$ containing $X$. The following statements are equivalent:
\begin{enumerate}
\item \label{1} $\mc{J}_{X/C}$ is Cohen-Macaulay,
\item \label{3} $\projdim(\derkxc)\leqslant k-1$,
\item \label{4} $\projdim(\derkxc)=k-1$.
\end{enumerate}

If in addition $X$ is Cohen-Macaulay then the previous properties are also equivalent to:
\begin{enumerate}
\setcounter{enumi}{3} 
\item \label{2}$\co_X/\mc{J}_{X/C}$ is Cohen-Macaulay of dimension $m-k-1$ or $\co_X/\mc{J}_{X/C}=(0)$.
\end{enumerate} 

\end{proposition}

Some equivalences are mentionned in \cite{gsci} for reduced complete intersections. In the case of a complete intersection $C$, the last characterization shows that freeness has a geometric interpretation since $C$ is free if and only if $C$ is smooth or the singular locus of $C$ is Cohen-Macaulay of codimension $1$ in $C$.

\begin{definition}\label{freeci}
A reduced equidimensional subspace $X$  is called \emph{free} if one of the equivalent properties \eqref{1}, \eqref{3}, \eqref{4} of proposition~\ref{char:gs} is satisfied.
\end{definition}

\smallskip

As it is mentioned in the introduction, it is also very interesting to consider the module of multi-logarithmic forms, which leads us to our main theorem~\ref{char:free}. We also deduce from theorem~\ref{char:free} a characterization of freeness involving logarithmic multi-residues (see corollary~\ref{nonci:cor:rc:free}).

\begin{theorem}
\label{char:free}
We keep the same hypothesis as in proposition~\ref{char:gs}. The following statements are equivalent:
 
\begin{enumerate}
\item \label{A} $X$ is free,
\item \label{B} $\projdim \lrp{\omkxc}\leqslant k-1$,
\item \label{C} $\projdim \lrp{\omkxc} = k-1$.
\end{enumerate} 
\end{theorem}

In particular, for $k=1$, we recognize the several characterizations of freeness for divisors we mentioned before.  In the hypersurface case, the duality between $\mathrm{Der}(-\log D)$ and $\Omega^1(\log D)$ gives immediately the fact that if one of the two modules is free, the other one is also free, whereas for equidimensional subspaces of codimension greater than $2$, the statement on the projective dimension of $\Omega^k(\log X/C)$ needs much more work.

\subsection{Proof of proposition~\ref{char:gs}}

The proof of the equivalences of proposition~\ref{char:gs} is based on the depth lemma as stated in \cite[Lemma 6.5.18]{dejong} and the Auslander-Buchsbaum formula.

Let us prove \eqref{1} $\Rightarrow$ \eqref{4}. We recall the following exact sequence: 
\begin{equation}
\label{derks}
 0\to \derkxc\to\Theta^k_S\to \mc{J}_{X/C}\to 0.
 \end{equation}

Then, thanks to the depth lemma, since the depth of $\mc{J}_{X/C}$ is $m-k$ and the depth of $\Theta^k_S$ is $m$, we have $\depth (\derkxc)=m-k+1$. By the Auslander-Buchsbaum formula, we have $\projdim (\derkxc)=k-1$.

\sm

The implication \eqref{4} $\Rightarrow$ \eqref{3} is trivial. 

\sm

Let us prove \eqref{3} $\Rightarrow$ \eqref{1}. By the Auslander-Buchsbaum formula, $\depth(\derkxc)\geqslant m-k+1$. In addition, we have $\depth(\Theta^k_S)=m$, and $\depth(\mc{J}_{X/C})\leqslant m-k$. As a consequence of the exact sequence~\eqref{derks} and of the depth lemma we have $\depth(\derkxc)=m-k+1$ and $\depth(\mc{J}_{X/C})=m-k$, so that $\mc{J}_{X/C}$ is maximal Cohen-Macaulay.

The equivalence between \eqref{1} and \eqref{2} is proved in \cite[Proposition 5.6]{snc} for Gorenstein spaces. Our proof for Cohen-Macaulay subspaces is completely similar.

If $\mc{J}_{X/C}=\co_X$, the statement is clear. Let us assume that $\mc{J}_{X/C}\neq \co_X$. 

Let us consider the following exact sequence of $\co_X$-modules: 
\begin{equation}
\label{jc}
0\to \mc{J}_{X/C}\to\mc{O}_X\to\co_X/\mc{J}_{X/C}\to 0.
\end{equation}
By assumption, $\co_X$ is Cohen-Macaulay of dimension $n$. Moreover, since we assume $C$ to be reduced, the singular locus of $C$ is of dimension at most $n-1$, and therefore the depth of $\co_X/\mc{J}_{X/C}$ is at most $n-1$. We deduce from the depth lemma that $\depth(\mc{J}_{X/C})=n \iff \depth(\co_X/\mc{J}_{X/C})=n-1$.

\sm 

\subsection{Preliminary to the proof of theorem~\ref{char:free}} 

Let us recall the following short exact sequence from proposition~\ref{nonci:prop:suite:ex}:
\begin{equation}
\label{free:eq:res}
0\to\tomf\to\omkxc\to\mc{R}_X\to 0.
\end{equation} 

The methods used to prove proposition~\ref{char:gs} applied to the short exact sequence~\eqref{free:eq:res} are not sufficient to prove directly theorem~\ref{char:free}.

The proof of theorem~\ref{char:free} is based on the explicit computation of some modules and morphisms of the long exact sequence obtained by applying the functor $\mathrm{Hom}_{\co_S}(-,\co_S)$ to the short exact sequence~\eqref{free:eq:res}: 
\begin{equation}
\label{ci:eq:long:seq}
0\to\mathrm{Hom}_{\co_S}(\mc{R}_X,\co_S)\to \mathrm{Hom}_{\co_S}(\omkxc,\co_S)\to \mathrm{Hom}_{\co_S}(\tomf,\co_S)\to\extt{1}{\mc{R}_X}{\co_S}\to\dots
\end{equation}

The structure of the proof is the following. Thanks to the Koszul complex, we compute the modules $\extt{q}{\totofq{k}}{\co_S}$. We then determine the modules $\extt{q}{\mc{R}_X}{\co_S}$ for $q\leqslant k$ using the change of rings spectral sequence. The most technical part is the explicit computation of the connecting morphism $$\alpha' :\extt{k-1}{\totofq{k}}{\co_S}\to \extt{k}{\mc{R}_X}{\co_S}.$$

This computation is necessary in order to identify the kernel and the image of $\alpha'$, which are used in the end of the proof. 

\subsubsection{}

We first compute the terms $\extt{q}{\totofq{k}}{\co_S}$ of the long exact sequence~\eqref{ci:eq:long:seq}. 

\begin{notation}
\label{nota:kos}
We denote by $K(\underline{f})$ the Koszul complex of $(f_1,\ldots,f_k)$ in $\co_S$: 
\begin{equation}
\label{suite:kos}
K(\underline{f}) \ :\  0\to\bigwedge^k\co_S^k\xrightarrow{d_k}\cdots\xrightarrow{d_{2}}\bigwedge^1\co_S^k\xrightarrow{d_1}\co_S\to 0.
\end{equation}

We also set $\wt{K}(\und{f})$ the complex obtained from $K(\und{f})$ by removing the last $\co_S$. 
\end{notation}

\begin{lemma}[\protect{\cite[Corollary 17.5, proposition 17.15]{eisenbudalgebra}}]
\label{kos:res}
Since the sequence $(f_1,\ldots,f_k)$ is regular, $K(\und{f})$ is a free $\co_S$-resolution of $\co_C$.

The dual complex $\mathrm{Hom}_{\co_S}(K(\und{f}),\co_S)$ of the Koszul complex is a free resolution of $\co_C$. 
\end{lemma}

\begin{remark}
\label{kos:res:sigma}
A consequence of lemma~\ref{kos:res} is that $\wt{K}(\und{f})$ gives a free $\co_S$-resolution of $\Sigma\simeq \mc{I}_C$.
\end{remark}

We can therefore use the complex $\wt{K}(\und{f})$ to compute the modules $\mathrm{Ext}^\bullet_{\co_S}\lrp{\tomf,\co_S}$:

\begin{lemma}
\label{lem:tom}We assume $k\geqslant 2$. 
The projective dimension of $\totofq{k}$ is $k-1$. Moreover, we have $\homss{\tomf}{\co_S}=f\Theta^k_S$, $\extt{k-1}{\tomf}{\co_S}=\Theta^k_S\otimes_{\co_S}\co_C$, and for all $j\notin\lra{0,k-1}$, $\extt{j}{\tomf}{\co_S}=0$.
\end{lemma}

\begin{proof} Since $\tomf=\Omega^k_S\otimes_{\co_S} \Sigma$, we have for all $q\in\N$, $\mathrm{Ext}^q_{\co_S}\lrp{\tomf,\co_S}= \Theta^k_S \otimes_{\co_S} \mathrm{Ext}^q_{\co_S}\lrp{\Sigma,\co_S}$. 

By remark~\ref{kos:res:sigma}, $\projdim(\tomf)=k-1$ and  $\extt{k-1}{\tomf}{\co_S}=\Theta^k_S\otimes_{\co_S}\co_C$ and for all $j\notin\lra{0,k-1}$, $\extt{j}{\tomf}{\co_S}=0$.

In addition, since $\extt{1}{\co_C}{\co_S}=0$, we deduce from the short exact sequence $$0\to\Sigma \xrightarrow{f\cdot}\co_S\to\co_C\to 0$$ that $\mathrm{Hom}_{\co_S}\lrp{\co_S,\co_S}$ and $\mathrm{Hom}_{\co_S}\lrp{\Sigma,\co_S}$ are isomorphic. More precisely, $\mathrm{Hom}_{\co_S}\lrp{\Sigma,\co_S}=f\co_S$. Hence the result.
\end{proof}

\subsubsection{}

We compute the modules $\extt{q}{\mc{R}_X}{\co_S}$ for $q\leqslant k$. 

To compute the modules involving $\mc{R}_X$, we introduce the change of rings spectral sequence (see for example \cite[Chapter XV and XVI]{cartaneilenberg} for details on spectral sequences). The change of rings spectral sequence applied to  an $\co_C$-module $M$ and $\co_S$ gives:
\begin{equation}
\label{suite:spec}
E_2^{pq}=\mathrm{Ext}_{\co_C}^{p}\lrp{M,\extt{q}{\co_C}{\co_S}}\Rightarrow \extt{p+q}{M}{\co_S}.
\end{equation}

\begin{lemma}
\label{spec:res}
For all $q<k$, $\extt{q}{M}{\co_S}=0$ and $\extt{k}{M}{\co_S}=\mathrm{Hom}_{\co_C}\lrp{M,\co_C}$. 
\end{lemma}
\begin{proof}
Since $(f_1,\ldots,f_k)$ is  a regular sequence, we have for all $q\neq k$, $\extt{q}{\co_C}{\co_S}=0$ and $\extt{k}{\co_C}{\co_S}=\co_C$. Therefore, the only non zero terms of the second sheet of the spectral sequence~\eqref{suite:spec} are the $E_2^{pk}$, so that the spectral sequence degenerates at rank $2$.  Hence the result. 
\end{proof}

We deduce from the previous propositions the following exact sequence: 
\begin{corollary}
\label{free:cor:seq}
The long exact sequence~\eqref{ci:eq:long:seq} gives: 
{\small
\begin{equation}
\label{free:eq:long:seq}
\dots\to 0\to \extt{k-1}{\omkxc}{\co_S}\to\Theta^k_S\otimes_{\co_S}\co_C\xrightarrow{\alpha} \mc{R}_X^\vee\to \extt{k}{\omkxc}{\co_S}\to 0\to \dots
\end{equation}}
where $\mc{R}_X^\vee=\mathrm{Hom}_{\co_C}\lrp{\mc{R}_X,\co_C}$. 
\end{corollary}

\subsubsection{Computation of the connecting morphism}
The previous results show that there exist isomorphisms $\beta$ and $\beta'$ such that the following diagram is commutative: 
\begin{center}
\begin{tikzpicture}[scale=0.3][description/.style={fill=white,inner sep=2pt}]
\matrix (m) [matrix of math nodes, row sep=1.5em,
column sep=0.8em, text height=2ex, text depth=0.25ex]
{ \Theta^k_S\otimes_{\co_S}\co_C &  & \mc{R}_X^\vee \\
\mathrm{Ext}^{k-1}_{\co_S}\lrp{\tomf,\co_S} & & \mathrm{Ext}^{k}_{\co_S}\lrp{\mc{R}_X,\co_S}\\};
\path[->,font=\scriptsize]
(m-2-1) \flecheb {$\beta'$} (m-1-1)
(m-2-3) \flecheb {$\beta$} (m-1-3)
(m-1-1) \fleche {$\alpha$} (m-1-3)
(m-2-1) \flecheb {$\alpha'$} (m-2-3);
\end{tikzpicture}
\end{center}

We recall that $c_X$ is the fundamental form of $X$ (see notation~\ref{nonci:nota:alpha0}). In particular, if $X$ is a complete intersection defined by $(h_1,\ldots,h_k)$, we have $c_X=\dd h_1\wedge\dots\wedge\dd h_k$.

The purpose of this subsection is to prove the following proposition: 
\begin{proposition}
\label{free:prop:alpha}
The connecting morphism of the exact sequence of corollary~\ref{free:cor:seq} is:
$$\begin{aligned}\alpha : 
\Theta^k_S\otimes_{\co_S}\co_C &\to\mc{R}_X^\vee\\
\delta\otimes \ov{a} & \mapsto a\cdot \delta(c_X)
\end{aligned}$$

In particular, the image of $\alpha$ is $\mc{J}_{X/C}$. 
\end{proposition}

Thanks to this proposition, we are able to compare $\mc{J}_{X/C}$ and $\mc{R}_X^\vee$, which is used in the end of the proof of theorem~\ref{char:free}. 

\sm

The computation of $\alpha$ is quite technical. We determine explicitly the isomorphisms $\beta$ and $\beta'$, and the connecting morphism $\alpha'$. 

\sm

We fix  an injective resolution $(\mc{I}^\bullet,\eps_\bullet)$ of $\co_S$. 

\begin{lemma}
\label{free:lem:beta}
Let $M$ be a finite type $\co_C$-module. The isomorphism of lemma~\ref{spec:res} is: 
{\small
\begin{align*}
\beta : \mathrm{Ext}^{k}_{\co_S}\lrp{M,\co_S}=H^k\lrp{\homss{M}{\mc{I}^\poin}} &\to \mathrm{Hom}_{\co_C}\lrp{M,H^k\lrp{\homss{\co_C}{\mc{I}^\poin}}}=\Hom{\co_C}{M}{\co_C}\\
[\psi]&\mapsto \lrp{ \wt{\psi} : \rho\mapsto [\wt{\psi}_\rho : \ov{a}\mapsto a.\psi(\rho)]}
\end{align*} }
\end{lemma}
\begin{proof}
Let $(P_p,\delta)$ be a free $\co_C$-resolution of $M$. There are two spectral sequences associated with the double complex $A^{pq}=\mathrm{Hom}_{\co_C}\lrp{P_p,\homss{\co_C}{\mc{I}^q}}$. The announced isomorphism follows from the definitions of the spectral sequences (see \cite[Chapter XV and XVI]{cartaneilenberg}) and the fact that both degenerate at rank two. 
\end{proof}

\begin{lemma}
\label{free:lem:tom:map}
The following map is the isomorphism of lemma~\ref{lem:tom}:
$$ \begin{aligned} \beta' : \underbrace{H^{k-1}\lrp{\Hom{\co_S}{\tomf}{\mc{I}^\bullet}}}_{=\extt{k-1}{\tomf}{\co_S}}&\to \Theta^k_S\otimes_{\co_S} \underbrace{H^{k-1}\lrp{\mc{I}^\bullet/\mathrm{Ann}_{\mc{I}^\bullet} (f_1,\ldots,f_k)}}_{=\co_C}\\
[\varphi] &\mapsto \sum_I \dr x_I\otimes [\ov{m_I}]
\end{aligned}$$
where $m_I\in\mc{I}^{k-1}$ satisfies $f\cdot m_I=\varphi(\dd x_I)$.
\end{lemma}
\begin{proof}
For all $j\in \N$, there is an isomorphism $\zeta:\Hom{\co_S}{\tomf}{\mc{I}^j}\to \Theta^k_S\otimes_{\co_S} \Hom{\co_S}{\Sigma}{\mc{I}^j}$ given by $\zeta(\varphi)=\sum_I \dr x_I \otimes \lrp{a\mapsto \varphi(a\dd x_I)}$.

Since $0\to \Sigma\xrightarrow{f} \co_S\to\co_C\to 0$ is exact and $\mc{I}^j$ is injective, the following map is an isomorphism: \begin{align*}
\homss{\co_S}{\mc{I}^j}/\homss{\co_C}{\mc{I}^j}&\to \homss{\Sigma}{\mc{I}^j}\\
  \left[\varphi : \co_S \to \mc{I}^j\right] &\mapsto  \left( a\mapsto \varphi(f\cdot a)\right)
\end{align*}

 Moreover, $\homss{\co_S}{\mc{I}^j}\simeq \mc{I}^j$ and $\homss{\co_C}{\mc{I}^j}\simeq \mathrm{Ann}_{\mc{I}^j}(f_1,\ldots, f_k)$, so that 
we obtain the isomorphism \begin{align*}
\xi : \mc{I}^j/\mathrm{Ann}_{\mc{I}^j}(f_1,\ldots,f_k)&\to \homss{\Sigma}{\mc{I}^j}\\
 [\ov{m}] &\mapsto  \left( a\mapsto a\cdot fm\right)
\end{align*}

Using the isomorphisms $\zeta$ and $\xi^{-1}$ we obtain the isomorphism $\beta'$ announced in the statement of lemma~\ref{free:lem:tom:map}.
\end{proof}

As we mention before in lemma~\ref{lem:tom}, $H^{k-1}\lrp{\mathrm{Hom}_{\co_S}\lrp{\Sigma,\co_S}}=\co_C$. Therefore, there exists an isomorphism 
$\gamma_1 : H^{k-1}\lrp{\mc{I}^\bullet/\mathrm{Ann}_{\mc{I}^\bullet}(f_1,\ldots,f_k)}\to \co_C$.
%}

Moreover, since for all $j\in\N$, $\mathrm{Ann}_{\mc{I}^j}(f_1,\ldots,f_k)$ is isomorphic to $\mathrm{Hom}_{\co_S}(\co_C,\mc{I}^j)$, we obtain an isomorphism 
$\gamma_2 : H^k\lrp{\mathrm{Ann}_{\mc{I}^\bullet} (f_1,\ldots,f_k)}\to \co_C$.

The following lemma gives the isomorphism between the modules $H^{k-1}\lrp{\mc{I}^\bullet/\mathrm{Ann}_{\mc{I}^\bullet}(f_1,\ldots,f_k)}$ and $H^k\lrp{\mathrm{Ann}_{\mc{I}^\bullet}(f_1,\ldots,f_k)}$.
\begin{lemma}
\label{free:lem:gamma}
The following map is an isomorphism:
\begin{align*}
\gamma : H^{k-1}\lrp{\mc{I}^\poin/\mathrm{Ann}_{\mc{I}^\poin} (f_1,\ldots,f_k)}&\to H^k\lrp{\mathrm{Ann}_{\mc{I}^\poin}(f_1,\ldots,f_k)}\\
 [\ov{m}]&\mapsto [\eps_{k-1}(m)]
\end{align*}
\end{lemma}
\begin{proof}
Let us denote  $\ov{\eps_{k-1}}:\mc{I}^{k-1}/\mathrm{Ann}_{\mc{I}^{k-1}} (f_1,\ldots,f_k)\to \mc{I}^k/\mathrm{Ann}_{\mc{I}^k} (f_1,\ldots,f_k)$. We first prove that $\gamma$ is well defined.

If $\ov{m}\in\mathrm{Ker}(\ov{\eps_{k-1}})$ then $\eps_{k-1}(m)\in\mathrm{Ann}_{\mc{I}^k} (f_1,\ldots,f_k)$. If $\ov{m}=\ov{\eps_{k-2}}(\ov{m'})$ for an element $\ov{m'}\in\mc{I}^{k-2}/\mathrm{Ann}_{\mc{I}^{k-2}} (f_1,\ldots,f_k)$, then $[\eps_{k-1}(\eps_{k-2}(m'))]=0$ so that the map $\gamma$ is well defined. 

Let us assume that $[\eps_{k-1}(m)]=0$. Then, there exists $m'\in\mathrm{Ann}_{\mc{I}^{k-1}} (f_1,\ldots,f_k)$ such that $\eps_{k-1}(m)=\eps_{k-1}(m')$, so that $m-m'\in\mathrm{Ker}(\eps_{k-1})=\mathrm{Im}(\eps_{k-2})$. Hence $[\ov{m}]=0$, therefore, the map $\gamma$ is injective. Let us consider $[m]\in H^k\lrp{\mathrm{Ann}_{\mc{I}^\poin}(f_1,\ldots,f_k)}$. Then $\eps_k(m)=0$ thus there exists $m'\in\mc{I}^{k-1}$ such that $\eps_{k-1}(m')=m$. Then $[m]=\gamma([\ov{m'}])$ so that $\gamma$ is surjective. 
\end{proof}

We now have all the identifications we need to compute $\alpha$. 

\begin{proof}[Proof of proposition~\ref{free:prop:alpha}]
Let us construct explicitly the connecting morphism:
$$\alpha' : H^{k-1}\lrp{\homss{\tomf}{\mc{I}^\poin}}\to H^k\lrp{\homss{\mc{R}_X}{\mc{I}^\poin}}.$$

We use a diagram chasing process based on the following commutative diagram: 
\begin{center}
\begin{tikzpicture}[scale=0.3][description/.style={fill=white,inner sep=2pt}]
\matrix (m) [matrix of math nodes, row sep=1.5em,
column sep=0.8em, text height=2ex, text depth=0.25ex]
{ 0 & & \homss{\tomf}{\mc{I}^{k-1}} & & \homss{\omkxc}{\mc{I}^{k-1}} & & \homss{\mc{R}_X}{\mc{I}^{k-1}} & & 0   \\
  0 & & \homss{\tomf}{\mc{I}^{k}} & & \homss{\omkxc}{\mc{I}^{k}} & & \homss{\mc{R}_X}{\mc{I}^{k}}  & & 0 \\};
\path[->,font=\scriptsize]
(m-1-3) \flecheb {} (m-1-1)
(m-1-5) \flecheb {$i^*$} (m-1-3)
(m-1-7) \flecheb {$\resxc^*$} (m-1-5)
(m-1-9) \flecheb {} (m-1-7)
(m-2-3) \flecheb {} (m-2-1)
(m-2-5) \flecheb {$i^*$} (m-2-3)
(m-2-7) \flecheb {$\resxc^*$} (m-2-5)
(m-2-9) \flecheb {} (m-2-7)
(m-1-3) \fleche {$\eps_{k-1}$} (m-2-3)
(m-1-5) \fleche {$\eps_{k-1}$} (m-2-5)
(m-1-7) \fleche {$\eps_{k-1}$} (m-2-7);
\end{tikzpicture}
\end{center}
Let $\varphi : \tomf\to\mc{I}^{k-1}$ be such that $\eps_{k-1}(\varphi)=0$. Let $\delta\otimes [\ov{m}]\in\Theta^k_S\otimes H^{k-1}\lrp{\mc{I}^\poin/\mathrm{Ann}_{\mc{I}^\poin}(f_1,\ldots,f_k)}$ be the image of $[\varphi]\in \mathrm{Ext}^{k-1}_{\co_S}\lrp{\tomf,\co_S}$ by $\beta'$. In particular, it means that for $\eta\in\tomf$, $\varphi(\eta)=\delta(f\eta)\cdot m$.

There exists $\Phi : \omkxc\to \mc{I}^{k-1}$ such that $\Phi\circ i=\varphi$. Let $\w\in\omkxc$. By proposition~\ref{nonci:prop:carac:loga}, there exists $g,\xi,\eta$ such that $g\w=\xi \frac{c_X}{f} +\eta$. Then $$g\Phi(\w)=\xi \Phi\lrp{\frac{c_X}{f}}+\varphi(\eta).$$

 Moreover,  for all $i\in\unk$, $f_i\Phi\lrp{\frac{c_X}{f}}=\varphi\lrp{f_i \frac{c_X}{f}}=f_i\delta(c_X)\cdot m$. Therefore, $$\Phi\lrp{\frac{c_X}{f}}=\delta(c_X)\cdot m+m'$$ with $m'\in\mathrm{Ann}_{\mc{I}^{k-1}}(f_1,\ldots, f_k)$. 

The image by $\eps_{k-1}$ of $\Phi$ satisfies:
$$g\cdot\eps_{k-1}(\Phi)(\w)=\xi \lrp{ \delta(c_X)\cdot \eps_{k-1}(m)+\eps_{k-1}(m')}.$$

Since $i^*\eps_{k-1}(\Phi)=0$, there exists $\Psi : \mc{R}_X\to \mc{I}^k$ such that $\eps_{k-1}(\Phi)=\resxc^*(\Psi)$. In particular, for all $\rho\in\mc{R}_X$, we have\footnote{We notice that $\eps_{k-1}(m)$ and $\eps_{k-1}(m')$ are canceled by $(f_1,\ldots,f_k)$, so that multiplying by $g\rho\in\co_C$ makes sense.}:
$$g\Psi(\rho)=g\rho \lrp{ \delta(\dd h_1\wedge\dots\wedge \dd h_k)\cdot \eps_{k-1}(m)+\eps_{k-1}(m')}.$$ 

We thus obtain the expression of $g\alpha'(\beta'^{-1}(\delta\otimes[\ov{m}]))\in\mathrm{Ext}^k_{\co_S}(\mc{R}_X,\co_S)$. 

By the isomorphism $\beta$ of lemma~\ref{free:lem:beta}, and using the identification of $\mathrm{Hom}_{\co_S}\lrp{\co_C,\mc{I}^\bullet}$ with $\mathrm{Ann}_{\mc{I}^\bullet}(f_1,\ldots,f_k)$, the class of $[g\Psi]\in H^k\lrp{\mathrm{Hom}_{\co_S}\lrp{\mc{R}_X,\mc{I}^\bullet}}$ corresponds to the map:
 \begin{align*}
\mc{R}_X &\to H^k\lrp{\mathrm{Ann}_{\mc{I}^\bullet}(f_1,\ldots,f_k)}\\
\rho&\mapsto [g\rho \cdot \lrp{ \delta(c_X)\cdot \eps_{k-1}(m)+\eps_{k-1}(m')}]
\end{align*} 

In addition, since $m'\in\mathrm{Ann}_{\mc{I}^{k-1}}(f_1,\ldots,f_k)$, we have for all $\rho\in\mc{R}_X$:
$$[g\rho \cdot \lrp{ \delta(c_X)\cdot \eps_{k-1}(m)+\eps_{k-1}(m')}]=[g\rho \cdot \lrp{ \delta(c_X)\cdot \eps_{k-1}(m)}]$$
 
Moreover, we have the isomorphisms: 

$$\co_C\xleftarrow{\gamma_1}  H^{k-1}\lrp{\mc{I}^\poin/\mathrm{Ann}_{\mc{I}^\poin} (f_1,\ldots,f_k)}\xrightarrow{\gamma} H^k\lrp{\mathrm{Ann}_{\mc{I}^\poin}(f_1,\ldots,f_k)}\xrightarrow{\gamma_2}\co_C$$

Let $\ov{a}=\gamma_1([\ov{m}])\in\co_C$. Since $\gamma, \gamma_1, \gamma_2$ are isomorphisms, we can assume that $\gamma_2\circ\gamma\circ\gamma_1^{-1}(\ov{1})=\ov{1}$, so that $\gamma_2([\eps_{k-1}(m)])=\ov{a}\in\co_C$. 

Consequently, $[g\Psi]$ is identified with the map $\left\{ \begin{aligned} 
\mc{R}_X& \to \co_C\\
\rho& \mapsto \ov{\rho g\delta(c_X) a }
\end{aligned}\right.$, and since $g$ is a non zero divisor in $\co_C$, the map $[\Psi]$ is identified with:
 $$\left. \begin{aligned} 
\mc{R}_X& \to \co_C\\
\rho& \mapsto \ov{\rho\delta(c_X) a }
\end{aligned}\right.$$

%}
Hence the result: let $\delta\otimes\ov{a}\in\Theta_S\otimes_{\co_S} \co_C$, then $\alpha(\delta\otimes \ov{a})=\ov{a}\cdot\delta(c_X)$.
\end{proof}

\subsection{End of the proof of theorem~\ref{char:free}}

We also need the following results, which are obtained from the following short exact sequence by using similar methods as the ones used in the proof of proposition~\ref{char:gs}:
$$0\to\tomf\to\omkxc\to\mc{R}_X\to 0.$$

We first notice the following property, which is a direct consequence of \cite[Theorem 21.21]{eisenbudalgebra}:
\begin{lemma}
\label{res:mcm}
If $X$ is a free equidimensional subspace, then $\mc{R}_X$ is a maximal Cohen-Macaulay module and $\mc{R}_X^\vee\simeq\mc{J}_{X/C}$. 
\end{lemma}

\begin{lemma}
\label{faible}
Let $X$ be a reduced equidimensional subspace. If $\projdim(\omkxc)\leqslant k-1$, then $\mc{R}_X$ is a maximal Cohen-Macaulay module. 
If $X$ is free, then $\projdim(\omkxc)\leqslant k$.
\end{lemma}

\begin{proof}
Let us consider the exact sequence $0\to \tomf\to\omkxc\to \mc{R}_X\to 0$. 

If $\projdim(\omkxc)\leqslant k-1$, by Auslander-Buchsbaum formula, $\depth(\omkxc)\geqslant m-k+1$. Since $\depth(\tomf)=m-k+1$ and $\depth(\mc{R}_X)\leqslant m-k$, by the depth lemma, $\depth(\mc{R}_X)=m-k=\dim(\mc{R}_X)$. 

If $X$ is free, by lemma~\ref{res:mcm}, we have $\depth(\mc{R}_X)=m-k$. By lemma~\ref{lem:tom} and Auslander Buchsbaum Formula, we have $\depth(\tomf)=m-k+1$. Therefore, by the depth lemma, $\depth(\omkxc)\geqslant m-k$ and $\projdim(\omkxc)\leqslant k$.
\end{proof}

Thanks to the explicit computation of the connecting morphism $\alpha$ of proposition~\ref{free:prop:alpha}, we are able to compare $\mathrm{Im}(\alpha)=\mc{J}_{X/C}$ and $\mc{R}_X^\vee$, so that we can finish the proof of theorem~\ref{char:free}, using lemma~\ref{faible}.

 \begin{proof}[End of the proof of theorem~\ref{char:free}]
 We start with the implication \eqref{A} $\Rightarrow$ \eqref{C}. By lemma~\ref{faible}, we have $\projdim(\omkxc)\leqslant k$. Moreover, by lemma~\ref{res:mcm}, $\mc{R}_X^\vee=\mc{J}_{X/C}$ so that the map $\alpha$ of proposition~\ref{free:prop:alpha} is surjective. Therefore, we have $\extt{k}{\omkxc}{\co_S}=0$.

Let $(\co_S^{\ell_j}, d_j)_{0\leqslant j\leqslant k}$ be a minimal free resolution of $\omkxc$. In particular, all the coefficients of $d_k$ belongs to the maximal ideal $\m$ of $\co_S$. The module  $\extt{k}{\omkxc}{\co_S}$ is isomorphic to $\co_S^{\ell_k}/\text{Im}({}^td_k)$, and is equal to zero.  By Nakayama lemma, $\co_S^{\ell_k}=0$ and therefore $\projdim(\omkxc)\leqslant k-1$.

In addition, since there are relations between the maximal minors of the Jacobian matrix, the map $\alpha$ has a non zero kernel.  Therefore, $\extt{k-1}{\omkxc}{\co_S}\neq 0$ and $\projdim(\omkxc)= k-1$.

\sm

The implication \eqref{C} $\Rightarrow$ \eqref{B} is trivial.

\sm

Let us prove \eqref{B} $\Rightarrow$ \eqref{A}. 

\sm

We assume $\mathrm{projdim} \lrp{\omkxc}\leqslant k-1$. The exact sequence \eqref{free:eq:long:seq} becomes:
$$0\to \extt{k-1}{\omkxc}{\co_S}\to\Theta^k_S\otimes_{\co_S}\co_C\xrightarrow{\alpha} \mc{R}_X^\vee\to 0.$$

Since by proposition~\ref{free:prop:alpha} the image of $\alpha$ is $\mc{J}_{X/C}$,  we have $\mc{R}_X^\vee=\mc{J}_{X/C}$.  By lemma~\ref{faible}, $\mc{R}_X$ is a maximal Cohen-Macaulay $\co_C$-module. Therefore, by \cite[Theorem 21.21]{eisenbudalgebra}, $\mc{J}_{X/C}$ is also maximal Cohen-Macaulay.
 \end{proof}

The following corollary gives other characterizations of freeness using the module of multi-residues:
\begin{corollary}
\label{nonci:cor:rc:free}
Let $X$ be a reduced equidimensional subspace contained in a reduced complete intersection $C$ of the same dimension. The following statements are equivalent:
\begin{enumerate}
\item $X$ is free,
\item $\projdim(\mc{R}_X)\neq \projdim(\omkxc)$,
\item $\mc{R}_X$ is Cohen-Macaulay and $\mathrm{Hom}_{\co_C}(\mc{R}_X,\co_C)\simeq\mc{J}_{X/C}$.
\end{enumerate}
\end{corollary}
\begin{proof}
We start with $(1)\iff (2)$. We consider the exact sequence~\eqref{free:eq:res}. The depth of $\tomf$ is $m-k+1$. Since $\depth(\mc{R}_X)\leqslant m-k$, the depth lemma gives that $\depth(\omkxc)\neq\depth(\mc{R}_X)$ if and only if $\depth(\mc{R}_X)=m-k$ and $\depth(\omkxc)\geqslant m-k+1$.  This is equivalent to the fact that $X$ is free by theorem~\ref{char:free} and the Auslander-Buchsbaum formula. 

The implication $(1)\Rightarrow (3)$ is given by lemma~\ref{res:mcm}. Conversely, if $\mathrm{Hom}_{\co_C}(\mc{R}_X,\co_C)=\mc{J}_{X/C}$ and $\mc{R}_X$ is Cohen-Macaulay, then, by \cite[Theorem 21.21]{eisenbudalgebra}, $\mc{J}_{X/C}$ is Cohen-Macaulay so that $X$ is free. 
\end{proof}

\begin{remark}
The condition that $\mc{R}_X$ is Cohen-Macaulay may not be satisfied (see for example \cite[Example 5.6]{orlik-terao-milnor-fibers}). 
\end{remark}
\subsection{Consequences of freeness}

We deduce from our results a computation of the $\mathrm{Ext}$-modules of both $\derkxc$ and $\omkxc$, which gives a relation between them which is more intricate than the $\Sigma$-duality considered in proposition~\ref{nonci:prop:dualite:omkxc:derkxc}.

\begin{corollary}
Let $X$ be an equidimensional  subspace of codimension at least 2. Then:
$$\Hom{\co_S}{\Omega^k(\log X/C)}{\co_S}=f\Theta^k_S,$$
$$\extq{k-1}{\co_S}{\Omega^k(\log X/C)}{\co_S}=\frac{\derkxc}{\lrp{\sum_{i=1}^k f_i \Theta^k_S}},$$
and for all $1\leqslant q\leqslant k-2$, $ \extq{q}{\co_S}{\Omega^k(\log X/C)}{\co_S}=0$. If moreover $X$ is free, then for all $q\geqslant k$, we have $ \extq{q}{\co_S}{\Omega^k(\log X/C)}{\co_S}=0$.
\end{corollary}
\begin{proof}
The beginning of the long exact sequence~\eqref{ci:eq:long:seq} gives  $$\Hom{\co_S}{\omkxc}{\co_S}\simeq \Hom{\co_S}{\tomf}{\co_S}=f\Theta^k_S.$$

Moreover, by lemma~\ref{spec:res}, for all $q\leqslant k-1$, $\extq{q}{\co_S}{\mc{R}_X}{\co_S}=0$ and for all $1\leqslant q\leqslant k-2$, $\extq{q}{\co_S}{\tomf}{\co_S}=0$ so that for all $1\leqslant q\leqslant k-2$, $\extq{q}{\co_S}{\omkxc}{\co_S}=0$. Thanks to the long exact sequence~\eqref{free:eq:long:seq}, we see that $\extq{k-1}{\co_S}{\omkxc}{\co_S}$ is the kernel of the map $\alpha: \Theta^k_S\otimes_{\co_S} \co_C\to \mc{J}_{X/C}$ computed in proposition~\ref{free:prop:alpha}. It is easy to see that this kernel is $\frac{\derkxc}{\lrp{\sum_{i=1}^k f_i \Theta^k_S}}$. In addition, if $X$ is free, $\projdim(\omkxc)=k-1$ which implies that for all $q>k-1$, $\extq{q}{\co_S}{\omkxc}{\co_S}=0$.
\end{proof}
\begin{remark}
If $X$ is not free, the module $\mathrm{Ext}^{k}_{\co_S}\lrp{\omkxc,\co_S}$ is isomorphic to $\mc{R}_{X}^\vee/\mc{J}_{X/C}$. 
\end{remark}

\begin{proposition}
Let $X$ be a reduced equidimensional subspace of codimension at least two. Then:
$$\Hom{\co_S}{\derkxc}{\co_S}=\Omega^k_S,$$
$$\mathrm{Ext}^{k-1}_{\co_S}\lrp{\derkxc,\co_S}\simeq \mc{R}_X,$$
and for all $1\leqslant q\leqslant k-2$, $\extt{q}{\derkxc}{\co_S}=0$. If moreover $X$ is free, then for all $q\geqslant k$, we have $\extt{q}{\derkxc}{\co_S}=0$.
\end{proposition}
\begin{proof}
We apply the functor $\mathrm{Hom}_{\co_S}(-,\co_S)$ to the following exact sequence:
$$0\to\derkxc\to\Theta^k_S\to\mc{J}_{X/C}\to 0$$

Since $\Theta^k_S$ is free, for all $q\geqslant 1$, $\mathrm{Ext}^q_{\co_S}(\Theta^k_S,\co_S)=0$. In addition, by lemma~\ref{spec:res}, for all $q<k$, $\mathrm{Ext}^q_{\co_S}(\mc{J}_{X/C},\co_S)=0$ and $\mathrm{Ext}^k_{\co_S}(\mc{J}_{X/C},\co_S)\simeq \mathrm{Hom}_{\co_S}\lrp{\mc{J}_{X/C},\co_C}=\mc{R}_X$. It gives us the result. 

If $X$ is free, by proposition~\ref{char:gs}, $\projdim(\derkxc)=k-1$ so that for all $q\geqslant k$, we have $\mathrm{Ext}^q_{\co_S}\lrp{\derkxc,\co_S}=0$. 
\end{proof}

\section{Computation of logarithmic modules for some families of free singularities}

\label{sec:4}
\label{curves}

In the case of free hypersurfaces, a natural question  is to determine a basis of the module of logarithmic forms or vector fields. The question of finding a generating family of the module of multi-logarithmic $k$-forms or $k$-vector fields also arises for reduced equidimensional subspaces, in addition, a new problem appears which is to determine the Betti numbers of a minimal free resolution of the module. It is completely done here in two cases. We first prove that normal crossing singularities are free singularities by computing explicitely the module $\derkxc$. We then compute explicitly a free resolution of the module of multi-logarithmic forms for  reduced quasi-homogeneous complete intersection curves.

\subsection{Normal crossing singularities of any codimension are free}
\label{normal-crossing-sec}

Let us first prove the following decomposition of multi-logarithmic vector fields, which is analogous to the hypersurface case:

\begin{proposition}
\label{prop:decomp:derlog}
Let $X$ be a reduced equidimensional subspace of codimension $k$, with irreducible components $X_1,\ldots, X_s$.
 Then:
$$\mathrm{Der}^k(-\log X)=\bigcap_{i=1}^s \mathrm{Der}^k(-\log X_i).$$
\end{proposition}

\begin{proof}
Let $C$ be a reduced complete intersection of codimension $k$ defined by a regular sequence $(f_1,\ldots,f_k)$ and which contains $X$. All the subspaces we will consider in this proof will be seen as subspaces of $C$.

Let $\delta\in\bigcap_{i=1}^s \mathrm{Der}^k(-\log X_i)$. Then for all $i\in\lra{1,\ldots,s}$, $\delta({c_{X_i}})\in\mc{I}_{X_i}$, where $c_{X_i}$ is a fundamental form for $X_i$ such that $\ov{c_{X_i}}=\ov{\beta\dd f_1\wedge \dots\wedge \dd f_k}\in\Omega^k\otimes_{\co_S}\co_C$ with $\beta$ satisfying the hypothesis of \ref{beta}.

In addition, from the definition of the fundamental form, one can see that $c_{X}=\sum_{i=1}^s c_{X_i}$. Therefore, $\delta(c_{X})=\sum_{i=1}^s \delta(c_{X_i})$. Furthermore, for all $i\in\lra{1,\ldots,s}$, the restriction of $\delta(c_{X_i})$ to any component of $C$ different from $X_i$ is zero. Therefore, $\delta(c_{X_i})\in \mc{I}_C$, and $\delta\in\mathrm{Der}^k(-\log X)$.  

Conversely, let $\delta\in\mathrm{Der}^k(-\log X)$. Let $i\in\lra{1,\ldots,s}$. Let us prove that $\delta\in\mathrm{Der}^k(-\log X_i)$. We have $\delta(c_{X})\in \mc{I}_C\subseteq \mc{I}_{X_i}$. In addition, as we already mentioned, for all $j\neq i$, $\delta(c_{X_j})\in \mc{I}_{X_i}$.Therefore, $\delta(c_{X_i})=\delta(c_X)-\sum_{j\neq i} \delta(c_{X_j})\in\mc{I}_{X_i}$ and we have $\delta\in\mathrm{Der}^k(-\log X_i)$. 
\end{proof}

Let us apply the previous proposition to the normal crossing singularities, which are generalization of the normal crossing divisor.

\begin{definition}[\protect{\cite[Definition 3.2]{sncv1}}]
\label{de:normal:crossing}
A normal crossing singularity is an arbitrary equidimensional union of coordinate subspaces.
\end{definition}

\begin{notation}
Let $E\subseteq \lra{1,\ldots,m}$. We denote by $X_E$ the coordinate subspace defined by the set of equations $\lra{ x_i=0\ ;\ i\in E}$. 
\end{notation}

\begin{proposition}
\label{prop:derk:normal:cross}
Let $k\in\lra{1,\ldots,m-1}$.  
Let $E_1,\ldots,E_r$ be pairwise distinct subsets of $\lra{1,\ldots,m}$ of cardinality $k$, with $E_j=\lra{i_1^j,\ldots, i_k^j}$. Let $X=\bigcup_{j=1}^r X_{E_j}$. Then:

$$\mathrm{Der}^k(-\log X)=\bigoplus_{j=1}^r \lrb{x_{i_1^j},\ldots,x_{i_k^j}}\dr_{E_j}\oplus \bigoplus_{\substack{|E|=k\\ E\notin\lra{E_1,\ldots,E_r}}} \co_S \dr_E.$$
\end{proposition}

Minimal free resolutions of $\lrb{x_{i_1^j},\ldots,x_{i_k^j}}$ are given by Koszul complex, so that we have that $\projdim\lrp{\lrb{x_{i_1^j},\ldots,x_{i_k^j}}}=k-1$ for all $j$. We thus obtain the following corollary, which shows in particular that the singular locus of the normal crossing divisor is free, which gives a positive answer to \cite[Question 1]{polthese}.
\begin{corollary}
Let $X$ be a normal crossing singularity. Then $X$ is a free singularity. 
\end{corollary}

\begin{proof}[Proof of proposition~\ref{prop:derk:normal:cross}]
We recall that the modules of multi-logarithmic vector fields are independent from the choice of the complete intersection containing the space. 

Let $j\in\lra{1,\ldots,r}$. The subspace $X_{E_j}$ is defined by the regular sequence $(x_i)_{i\in E_j}$, so that 
\begin{align*}
\mathrm{Der}(-\log X_{E_j})&=\lra{\delta\in\Theta^k\ ;\ \delta\lrp{\wedge_{i\in E_j} \dd x_i}\in \lrb{x_{i_1^j},\ldots,x_{i_k^j}}}\\
 &=\lra{\delta=\sum_{|E|=k} \delta_E\dr_E\in\Theta^k\ ;\ \delta_{E_j}\in\lrb{x_{i_1^j},\ldots,x_{i_k^j}}}
\end{align*}
Proposition~\ref{prop:decomp:derlog} then gives the result. 
\end{proof}

\subsection{Quasi-homogeneous curves}

Let us notice the following property  which is easy to prove from the definition of freeness: 

\begin{proposition}
Reduced curves in $(\C^m,0)$ are free  subspaces.
\end{proposition}

We use in this section our main theorem~\ref{char:free} and results from \cite{polcras} and \cite{polvalues}. In these papers, the author proves that the set of values of the module of multi-residues $\mc{R}_C$  satisfies a symmetry with the values of the Jacobian ideal, and gives the relation between the values of $\mc{R}_C$ and the values of the K\"ahler differentials for complete intersection curves. This result is then generalized in \cite{kts} for more general curves by considering the dualizing module.

\medskip

We describe here explicitly the module of multi-logarithmic differential forms for a quasi-homogeneous complete intersection curve.  

We recall the following notations from \cite{polvalues}.

\medskip

Let $C=C_1\cup \dots \cup C_p$ be a reduced complete intersection curve with $p$ irreducible components.  The normalization of $C$ satisfies $\co_{\wt{C}}=\bigoplus_{i=1}^p \C\lra{t_i}$. It induces for all $i\in\lra{1,\ldots,p}$ a valuation map $$\val_i : \mc{M}_C\ni g\mapsto \val_i(g)\in\Z\cup\lra{\infty}.$$

The value of an element $g\in\mc{M}_C$ is $\val(g)=(\val_1(g),\ldots,\val_p(g))\in \lrp{\Z\cup\lra{\infty}}^p$. For a fractional ideal $I\subset \mc{M}_C$, we set $\val(I):=\lra{\val(g)\ ;\ g\in I \text{ non zero divisor}}\subset \Z^p$.

We consider the product order on $\Z^p$, so that for all $\alpha,\beta\in\Z^p$, $\alpha\leqslant \beta$ means that for all $i\in\lra{1,\ldots,p}$, $\alpha_i\leqslant \beta_i$.  We set $\und{1}=(1,\ldots,1)$. 

We denote by $\mc{C}_C=\co_{\wt{C}}^\vee$ the conductor ideal. There exists $\gamma\in\N^p$ such that $\mc{C}_C=\und{t}^{\gamma}\co_{\wt{C}}$. In particular, $\gamma=\inf\lra{\alpha\in\N^p ; \alpha+\N^p\subseteq \val(\co_C)}$, and is called the conductor of $C$.

\medskip

Let $C$ be a reduced complete intersection curve defined by a regular sequence $(h_1,\ldots,h_{m-1})$. Let us consider the following properties:
\begin{condition}\ 
\label{cond:qh}
\begin{enumerate}[a)]
\item There exist $(w_1,\ldots,w_m)\in\N^m$ such that for all $i\in\lra{1,\ldots,m-1}$, $h_i$ is quasi-homogeneous of degree $d_i$ with respect to the weight $(w_1,\ldots,w_m)$. 
\item $m$ is the embedding dimension. Equivalently, for all $i\in\lra{1,\ldots,m-1}$, $h_i\in\m^2$ where $\m$ is the maximal ideal of $\co_S$.
\end{enumerate} 
\end{condition}

\subsubsection{Generators of $\Omega^{m-1}(\log C)$}

\begin{lemma}
\label{w0}
Let $C$ be a reduced complete intersection curve satisfying condition a), and $D$ be the hypersurface defined by $h$. Then $$\w_0=\frac{\sum_{i=1}^m (-1)^{i-1} w_i x_i \wh{\dd x_i}}{h}\in\Omega^{m-1}(\log D).$$
\end{lemma}
\begin{proof}
Let $i\in\lra{1,\ldots,m-1}$. We have $\sum_{k=1}^m w_kx_k \dfrac{\dr h_i}{\dr x_k}=d_i h_i$ so that $\dd h_i\wedge \w_0=\dfrac{d_i}{\wh{h_i}}\dd \und{x}$. Since $\dd h=\sum_{i=1}^{m-1} \wh{h_i} \dd h_i$, we have $\dd h\wedge \w_0\in\Omega^m_S$. Thus, $\w_0\in\Omega^{m-1}(\log D)$. 
\end{proof}

\begin{remark}
By remark~\ref{rem:incl:log}, we also have $\w\in\Omega^{m-1}(\log C)$. 
\end{remark}

For $i\in\lra{1,\ldots,m}$, we denote by $J_i$ the $(m-1)\times(m-1)$ minor of the Jacobian matrix obtained by removing the column $\lrp{\frac{\dr h_j}{\dr x_i}}_{1\leqslant j\leqslant m-1}$. 
\begin{lemma}
\label{res:lem:res:w0}
For all $i\in\lra{1,\ldots,m}$, $J_i\res(\w_0)=(-1)^{i-1}w_ix_i$. 
Let $c_1,\ldots,c_m\in\C$ be such that $\sum_{i=1}^m c_iJ_i$ induces a non zero divisor of $\co_C$. We thus  have $$\res(\w_0)=\frac{\sum_{i=1}^m (-1)^{i-1}c_i w_i x_i}{\sum_{i=1}^m c_i J_i}.$$
\end{lemma}
\begin{proof}
Let $i_0\in\lra{1,\ldots,m}$. We recall that $\dd h_1\wedge \dots \wedge \dd h_{m-1}=\sum_{i=1}^m J_i \wh{\dd x_i}$. We have:
\begin{equation}
\label{cur:eq:min}
J_{i_0} \w_0= (-1)^{i_0-1} w_{i_0} x_{i_0} \frac{\dd h_1\wedge \dots \wedge \dd h_{m-1}}{h}+\sum_{\substack{i=1 \\ i\neq i_0}}^m \underbrace{\lrp{(-1)^{i-1} w_i x_i J_{i_0} -(-1)^{i_0-1}w_{i_0}x_{i_0} J_i}}_{=:\lambda_i}\frac{\wh{\dd x_i}}{h}.
\end{equation}

Let  $i\in\lra{1,\ldots,m}$,  $i\neq i_0$. We develop $J_i$ with respect to the column $i_0$, and $J_{i_0}$ with respect to the column $i$.  For $\lra{i_1,i_2}\subseteq \lra{1,\ldots,m}$ and $j\in\lra{1,\ldots,m-1}$, we denote by $J_{i_1,i_2}^j$ the minor of the Jacobian matrix obtained by removing the columns $i_1,i_2$ and the line $j$. We then obtain:
$$\lambda_i=\mathrm{sgn}(i_0-i)\sum_{\ell =1}^{m-1} (-1)^{\ell-1} \lrp{ w_i x_i \frac{\dr h_\ell}{\dr x_i} + w_{i_0} x_{i_0} \frac{\dr h_\ell }{\dr x_{i_0}}}\cdot J_{i_0,i}^\ell$$

By developing a convenient determinant, one can prove that for $p\notin{i,i_0}$ and $\ell\in\lra{1,\ldots,m-1}$: 
\begin{equation}
\label{cur:eq:jij:zero}
\sum_{\ell=1}^{m-1} (-1)^{\ell-1}\frac{\dr h_\ell}{\dr x_p} J_{i,i_0}^\ell=0.
\end{equation}
We then have:
$$\lambda_i= \mathrm{sgn}(i_0-i) \sum_{\ell=1}^{m-1} (-1)^{\ell-1} \lrp{\sum_{p=1}^m w_p x_p \frac{\dr h_\ell}{\dr x_p}}J_{i_0,i}^\ell=\mathrm{sgn}(i_0-i) \sum_{\ell=1}^{m-1} (-1)^{\ell-1} d_\ell h_\ell J_{i_0,i}^\ell.$$

Therefore, there exists $\eta\in\toto{m-1}$ such that \begin{equation}
\label{res:eq:jio:wo}
J_{i_0}\w_0=(-1)^{i_0-1}w_{i_0}x_{i_0} \dfrac{\dd h_1\wedge\dots\wedge \dd h_k}{h_1\cdots h_k}+\eta.
\end{equation} 

 Hence the result. 
\end{proof}

\begin{lemma}
\label{val:w0}
We recall that $\gamma\in\N^p$ denotes the conductor of $C$.
With the notations of lemma~\ref{w0}, we have: $$\inf(\val(\mc{R}_C))=\val(\res(\w_0))=-\gamma+\und{1}$$
\end{lemma}
\begin{proof}
By \cite[Proposition 3.30 and (18)]{polvalues}, we have $\val(J_i)=\gamma+\val(x_i)-\und{1}$. By considering for all branch $C_i$ of $C$ an index $j(i)$ such that ${x_{i(j)}}|_{C_{i}}\neq 0$, we deduce from lemma~\ref{res:lem:res:w0} that $\val(\res(\w_0))=-\gamma+\und{1}$.  
  
Let us prove that $\inf(\val(\mc{R}_C))=\val(\res(\w_0))$. As in~\cite{polvalues}, we set for $v\in\Z^p$, $\Delta_i(v,\mc{J}_C)=\lra{\alpha\in\val(\mc{J}_C) ; \alpha_i=v_i \text{ and } \forall j\neq i, \alpha_j>v_j}$ and $\Delta(v,I)=\bigcup_{i=1}^p \Delta_i(v,I)$. 

By \cite[Proposition 3.30]{polvalues}, since $\val(\co_C)\subseteq \val(\Omega^1_C)$ (see \cite{hefez}), we also have 
 $\gamma+\val(\co_C)-\und{1}\subseteq \val(\mc{J}_C)$. Therefore, $2\gamma-\und{1}+\N^p\subseteq \val(\mc{J}_C)$.

We then have $\max\lra{v\in\Z^p ; \Delta(v,\mc{J}_C)=\emptyset}\leqslant 2\gamma-\und{2}$.

By \cite[Theorem 2.4]{polvalues} we have $v\in\val(\mc{R}_C)\iff \Delta(\gamma-v-\und{1},\mc{J}_C)=\emptyset$.
 It implies that $\inf(\val(\mc{R}_C))\geqslant -\gamma+\und{1}$. 
Hence the result: $\inf(\val(\mc{R}_C))=\val(\res(\w_0))$. 
\end{proof}

\begin{proposition}
\label{gen:rc}
Let $C$ be a singular complete intersection satisfying condition a). Then $\mc{R}_C$ is generated by $\res\lrp{\frac{\ddh}{h}}=1$ and $\res(\w_0)$, where $\w_0$ is given in lemma~\ref{w0}. In addition, this generating family is minimal. 
\end{proposition}
\begin{proof}
We set $Z=\text{Sing}(C)$ the singular locus of $C$. By dualizing over $\co_C$ the exact sequence $0\to\mc{J}_C\to\co_C\to \co_Z\to 0$, we obtain \begin{equation}
\label{KW}
0\to\co_C\to \mc{R}_C\to \w_Z\to 0
\end{equation}
 where $\w_Z$ is the dualizing module of $Z$. Moreover, the singular locus of a quasi-homogeneous curve is Gorenstein (see \cite[Satz 2]{kunzwaldijac}), so that $\w_Z=\co_Z$. The  exact sequence~\eqref{KW} implies that $\mc{R}_C$ is generated by two elements, the image of $1\in\co_C$, which is $1\in\mc{R}_C$, and the antecedent of $1\in\co_Z$. Therefore, there exists $\rho_0\in\mc{R}_C$ such that $(1,\rho_0)$ generates $\mc{R}_C$. 
 
  It remains to prove that we can take $\rho_0=\res(\w_0)$. By lemma~\ref{val:w0},  $\val(\res(\w_0))=-\gamma+\und{1}$. 
 
 We assume first that $-\gamma+\und{1}\notin \N^p$. For example, $-\gamma_1+1<0$.  There exists $\alpha_0,\alpha_1\in\co_C$ such that $\res(\w_0)=\alpha_0 \rho_0+\alpha_1$. Since $\val(\alpha_1)\geqslant 0$, and $\inf(\val(\mc{R}_C))=\val(\res(\w_0))$, we have $\val_1(\rho_0)=\val_1(\w_0)$ therefore $\val_1(\alpha_0)=0$ which implies that $\val(\alpha_0)=0$ and $\alpha_0$ is invertible. Thus, $(\res(\w_0),1)$ generates $\mc{R}_C$. 
 
 Let us assume that $-\gamma+\und{1}\in\N^p$. Since $\gamma\geqslant 0$, we must have $\gamma=\und{1}$ or $\gamma=\und{0}$. However, if $\gamma=\und{0}$, we have $\co_C=\co_{\wt{C}}$ so that $C$ is smooth. Therefore, $\gamma=\und{1}$. By \cite[Proposition 3.31]{polvalues},  we have $\val(\mc{J}_C)=\und{1}+\N^p=\val(\mc{C}_C)$.  It implies that $\mc{J}_C=\mc{C}_C$, so that by duality, $\mc{R}_C=\co_{\wt{C}}$. By \cite[Proposition 4.11]{snc}, it implies that $C$ is a plane normal crossing curve. By Saito criterion \cite[(1.8)]{saitolog}, if $h=xy$ defines a plane curve $C$, then $(\w_0=\frac{x\dd y-y\dd x}{h},\frac{\dd h}{h})$ is a basis of $\Omega^1(\log C)$. Hence the result. 

\end{proof}

Since $\lra{\frac{h_i\wh{\dd x_j}}{h} ; i\in\lra{1,\ldots,m-1},j\in\lra{1,\ldots,m}}$ generates $\wt{\Omega}^{m-1}_{\und{h}}$, proposition~\ref{gen:rc} gives:

\begin{corollary}
\label{gen}
Let $C$ be a singular complete intersection curve satisfying condition a). Then  $\Omega^{m-1}(\log C)$ is generated by the multi-logarithmic form $\w_0$ of lemma~\ref{w0},  $\dfrac{\ddh}{h}$, and the family $\lra{\frac{h_i\wh{\dd x_j}}{h} ; i\in\lra{1,\ldots,m-1},j\in\lra{1,\ldots,m}}$. 
\end{corollary}

We will see in the next section that this generating family is minimal if the conditions~\ref{cond:qh} are satisfied.

\subsubsection{Free resolution}
Since curves are free, lemma~\ref{res:mcm} yields $\depth(\mc{R}_C)=1$, so that by the Auslander-Buchsbaum Formula, the projective dimension of $\mc{R}_C$ as a $\co_S$-module is $m-1$. In addition, by theorem~\ref{char:free}, the projective dimension of $\Omega^{m-1}(\log C)$ is $m-2$. 

\sm

In corollary~\ref{gen}, we give a generating family of $\Omega^{m-1}(\log C)$. We can  go further and compute explicitly a free $\co_S$-resolution of $\mc{R}_C$ and $\Omega^{m-1}(\log C)$ for a quasi-homogeneous complete intersection curve. 

\begin{theorem}
\label{res:theo:res:res}
Let $C$ be a reduced complete intersection curve satisfying the conditions~\ref{cond:qh}. We set for $p\in\lra{0,\ldots,m-2}$: 
$$F_p=\bigwedge^p\co_{S}^{m-1}\oplus \bigwedge^p \co_S^m$$
and $F_{m-1}=\bigwedge^{m-1}\co_S^m$.

There exist differentials $\delta_\bullet : F_\bullet\to F_{\bullet-1}$ such that $(F_\bullet,\delta_\bullet)$ is a minimal free resolution of $\mc{R}_C$ as a $\co_S$-module.

In particular, the Betti numbers of $\mc{R}_C$ as a $\co_S$-module are:
$$\forall p\in\lra{0,\ldots,m-2}, b_j(\mc{R}_C)=\binom{m-1}{p}+\binom{m}{p}  \text{ and } b_{m-1}(\mc{R}_C)=m$$
\end{theorem}

In order to prove this theorem, we introduce the following exact sequence, where the middle module is isomorphic to  $\mc{R}_C$ and where a free resolution of the two other modules is given by Koszul complexes. 

\begin{lemma}
\label{res:lem:suite:ex}
We keep the hypothesis of theorem~\ref{res:theo:res:res}. Let $c_1,\ldots, c_m\in\C$ be such that $g=\sum_{i=1}^m c_i J_i\in\co_S$ induces a non zero divisor in $\co_C$. Let $y=\sum_{i=1}^m (-1)^{i-1}c_i w_i x_i$.  In particular, $\mc{R}_C\simeq y\co_C+g\co_C$. We have the following exact sequence:
\begin{equation}
\label{res:eq:rc:j1}
0\to \co_C\xrightarrow{y} y\co_C+g\co_C\to \frac{\lrp{g\co_S+y\co_S+\mc{I}_C}}{y\co_S+\mc{I}_C}\to 0
\end{equation}

In addition, a free resolution of $\dfrac{\lrp{g\co_S+y\co_S+ \mc{I}_C}}{y\co_S+\mc{I}_C}$ is given by the Koszul complex associated with the regular sequence $(w_1x_1,\ldots,w_mx_m)$. 
\end{lemma}

\begin{proof}
The exact sequence~\eqref{res:eq:rc:j1} is just a consequence of the fact that $y$ is a non zero divisor of $\co_C$.

We have:
$$0\to \lrp{\lrp{y,h_1,\ldots,h_{m-1}} : g}_{\co_S}\to \co_S\xrightarrow{g} \frac{\lrp{g\co_S+y\co_S+\mc{I}_C}}{y\co_S+\mc{I}_C}\to 0.$$

By~\eqref{res:eq:jio:wo}, for all $i\in\lra{1,\ldots,m}$, $J_i\w_0=(-1)^{i-1}w_ix_i\dfrac{\dd h_1\wedge \dots\wedge\dd h_{m-1}}{h}+\eta_i$ with $\eta_i\in\toto{m-1}$. We have for all $i,j\in\lra{1,\ldots,m}$:
\begin{equation}
\label{res:eq:ji:j1}
(-1)^{i-1}w_ix_iJ_j=(-1)^{j-1}w_jx_jJ_i \mod (h_1,\ldots,h_{m-1})
\end{equation}
We thus have $yJ_j=(-1)^{j-1}w_jx_j g \mod (h_1,\ldots,h_{m-1})$. Therefore, $$(w_1x_1,\ldots,w_mx_m)\subseteq \lrp{\lrp{y, h_1,\ldots,h_{m-1}} : g}_{\co_S}.$$ Moreover, by proposition~\ref{gen:rc}, $(1,\frac{y}{g})$ is a minimal generating family of $\mc{R}_C$, thus, $g\notin (y,h_1,\ldots,h_{m-1})$. We thus have
$\m=(w_1x_1,\ldots,w_mx_m)= \lrp{\lrp{y, h_1,\ldots,h_{m-2}} : J_1}_{\co_S}$. Hence the result.
\end{proof}

\begin{proof}[Proof of theorem~\ref{res:theo:res:res}]
We consider the exact sequence~\eqref{res:eq:rc:j1}. A minimal free resolution of $\co_C$ is given by the Koszul complex associated with the regular sequence $(h_1,\ldots,h_{m-1})$ and a minimal free resolution of $\frac{\lrp{g\co_S+y\co_S+\mc{I}_C}}{y\co_S+\mc{I}_C}$ is given by the Koszul complex associated to $(w_1x_1,\ldots,w_mx_m)$. We deduce from these two resolutions a free resolution $(F'_j,\delta'_j)$ of $(y,g)\co_C$, whose length is $m$. However, the projective dimension of $\mc{R}_C$ is $m-1$, so that the free resolution we obtain is not minimal. Thanks to the explicit computation of the differential $\delta_j'$, one can see that all the coefficients of the differentials belongs to the maximal ideal $\m$ of $\co_S$, except from $\delta_m'$ which has an invertible coefficient. A minimization of these free resolution then gives the announced result. The precise expression of the differential can be computed exactly as in \cite[Th\'eor\`eme 6.1.29]{polthese}, where we assume that $g=J_1$ is a non zero divisor in $\co_C$.
\end{proof}

\begin{theorem}
\label{res:theo:omk:free:res}
We keep the notations and hypothesis of theorem~\ref{res:theo:res:res}. We set for all $j\in\lra{0,\ldots,m-3}$, 
$$P_j=\lrp{\bigwedge^{j+1}\co_S^{m-1}\otimes\Omega_S^{m-1}}\oplus F_j$$
and $P_{m-2}=F_{m-2}$. There exist $\alpha_\bullet : P_\bullet\to P_{\bullet-1}$ such that $(P_\bullet,\alpha_\bullet)$ is a minimal free resolution of $\Omega^{m-1}(\log C)$.

In particular, the Betti numbers of $\Omega^{m-1}(\log C)$ are for all $j\in\lra{0,\ldots, m-3}$,  $$b_j(\Omega^{m-1}(\log C))=m\binom{m-1}{j+1}+\binom{m-1}{j}+\binom{m}{j} \text{ and } b_{m-2}(\Omega^{m-1}(\log C))=m-1+\binom{m}{m-2}.$$

\end{theorem}
\begin{proof}
We consider the exact sequence $$0\to\toto{m-1}\to \Omega^{m-1}(\log C)\to \mc{R}_C\to 0.$$

The free resolutions of $\toto{m-1}$ and $\mc{R}_C$ induce a free resolution $(P'_\bullet,\alpha'_\bullet)$ of $\Omega^{m-1}(\log C)$, whose length is $m-1$. By theorem~\ref{char:free}, the projective dimension of $\Omega^{m-1}(\log C)$ is $m-2$ since $C$ is free. Therefore, the previous free resolution of $\Omega^{m-1}(\log C)$ is not minimal. A minimization gives the announced result. The expression of the differentials can be found in \cite[Th\'eor\`eme 6.1.33]{polthese}
\end{proof}

\begin{remark}
In the statement of theorem~\ref{res:theo:omk:free:res}, we assume that $m$ is the embedded dimension. Let us assume that there exists $\ell\in\lra{1,\ldots,m-1}$ such that we have that for all $i\in\lra{1,\ldots,\ell}$, $h_i=x_i$ and for $i\geqslant \ell+1$, $h_i\in\C\{x_{\ell+1},\ldots,x_m\}\cap \m^2$. We denote by $D_i$ the hypersurface defined by $h_i$. 
We set $C'=D_{\ell+1}\cap\ldots\cap D_{m-1}\subseteq \C^{m-\ell}$, so that $C=\lra{0}\times C'$. Then one can prove that for all $q\geqslant \ell$ (see \cite[Proposition 3.1.24]{polthese}): 
$$\Omega^q(\log C)=\frac{1}{x_1\cdots x_\ell}\Omega^{q-\ell}(\log C')\wedge \dd x_1\wedge \dots \wedge \dd x_\ell+\frac{1}{h_1\cdots h_{m-1}}\mc{I}_C\Omega^q.$$
In particular, for all $q\geqslant k$, the module $\mc{R}_{C}^{q-k}$ is equal to the $\co_C$-module $\mc{R}_{C'}^{q-k}$. 
\end{remark}

\subsection{Examples and remarks}
\label{sec:ex}

We give in this subsection several examples and remarks on the properties of the modules of multi-logarithmic forms for more general equations and subspaces

 \begin{example}
\label{cur:ex:cm}
Let us consider the case of a Cohen-Macaulay quasi-homogeneous curve. Let us notice that for a quasi-homogeneous complete intersection curve in $(\C^3,0)$, a free resolution of $\Omega^2(\log C)$ is:
$$0\to \co_S^5\to\co_S^{8}\to\Omega^2(\log C)\to 0.$$
Let us consider the curve $X$ of $\C^3$ parametrized by $(t^3,t^4,t^5)$. This curve is a quasi-homogeneous curve which is not Gorenstein. The reduced ideal defining $X$ is the ideal generated by $h_1=xz-y^2$, $h_2=x^3-yz$ and $h_3=x^2y-z^2$. We set $C$ the reduced complete intersection defined by $(xz-y^2,x^3-yz)$. A computation made with \textsc{Singular} gives the following minimal free resolution of $\Omega^2(\log X/C)$ :
$$0\to \co_S^6\to \co_S^9\to \Omega^2(\log X/C)\to 0.$$

In particular, the number of generators of $\mc{R}_{X}$ is $3$.
\end{example}

\begin{remark}
If $C=C_1\cup \dots \cup C_p$ is a reduced quasi-homogeneous complete intersection curve, using theorem~\ref{res:theo:res:res}, we computed in \cite{pol-solomon-terao} a generating family of the module of multi-residues of any union $\bigcup_{i\in I} C_i$ with $I\subseteq\lra{1,\ldots,p}$.
\end{remark}

If the curve is not quasi-homogeneous, the number of generators can be strictly greater:

\begin{example}
Let $h_1=x_1^7-x_2^5+x_1^5x_2^3$ and $h_2=x_1^3x_2-x_3^2$. The sequence $(h_1,h_2)$ defines a reduced complete intersection curve of $\C^3$, which is not quasi-homogeneous. We use \textsc{Singular} to compute a minimal free resolution of $\Omega^{2}(\log C)$: 
$$0\to \co_S^6\to \co_S^9\to \Omega^{2}(\log C)\to 0.$$
\end{example}

\begin{remark}
Let $C$ be a reduced  complete intersection curve defined by a regular sequence $(h_1,\ldots,h_{m-1})$ such that the equations $h_1,\ldots,h_{m-1}$ are quasi-homogeneous for the same weights. By remark~\ref{rem:incl:log}, if $D$ denotes the hypersurface defined by $h=h_1\cdots h_k$, we have $\Omega^{m-1}(\log D)\subseteq \Omega^{m-1}(\log C)$. Proposition~\ref{gen:rc} shows that we have $\mathrm{res}_C(\Omega^{m-1}(\log D))=\mc{R}_C$. However, this property may not be satisfied for arbitrary equations of $C$. Indeed, let us assume for example that $\text{w-deg}(h_{1})\neq \text{w-deg}(h_{2})$. We set $f_1=h_1+h_2$ and for all $i\in\lra{2,\ldots,m-1}$, $f_i=h_i$. The ideal generated by $f_1,\ldots,f_{m-1}$ is $\mc{I}_C=(h_1,\ldots,h_{m-1})$. By lemma~\ref{w0} and remark~\ref{numerators}, $\w=\frac{\sum_{i=1}^m (-1)^{i-1} w_i x_i \wh{\dd x_i}}{f_1\ldots f_{m-1}}\in \Omega^{m-1}(\log C, \underline{f})$. It is possible to prove that $\res(\w)\notin\res(\Omega^{m-1}(\log D'))$, where $D'$ is the hypersurface defined by $f=f_1\cdots f_{m-1}$ (see \cite[Proposition 6.1.38]{polthese}). An interesting question would be to determine under which hypothesis the map $\mathrm{res}_C : \Omega^q(\log D)\to \mc{R}_C^{q-k}$ of \cite[\S 6, Theorem 2]{alekres} is indeed surjective.

\end{remark}

\begin{example}
\label{ex:surf:c4}
It is then natural to look for analogous results for quasi-homogeneous complete intersections of higher dimension. However, the situation seems to be more complicated than in the curve case, as it is shown with the following examples coming from subspace arrangements.

Let us consider the surface $C_1$ in $(\C^4,0)$ defined by the ideal $(xy,zt)$. In particular, $C_1$ is a reduced complete intersection. We compute with \textsc{Singular} a minimal free resolution of $\Omega^2(\log C_1)$:
$$0\to \co_S^{10}\to \co_S^{16}\to \Omega^2(\log C_1)\to 0.$$

For the reduced complete intersection surface $C_2$ of $(\C^4,0)$ defined by $(xy(x+y+z),zt)$, a minimal free resolution is:
$$0\to \co_S^{11}\to \co_S^{17}\to \Omega^2(\log C_2)\to 0.$$
\end{example}

\bibliographystyle{alpha-fr}
\bibliography{bibli2}
\end{document}